\documentclass[11pt]{article}
\usepackage{amssymb}
\usepackage{amsmath}
\usepackage{mathrsfs}
\usepackage{graphics}
\usepackage{graphicx}
\usepackage{xcolor}
\usepackage{subfigure}
\usepackage[T1]{fontenc}
\usepackage{latexsym,amssymb,amsmath,amsfonts,amsthm}\usepackage{txfonts}
\newcommand{\be}{\begin{eqnarray}}
\newcommand{\ben}{\begin{eqnarray*}}
\newcommand{\en}{\end{eqnarray}}
\newcommand{\enn}{\end{eqnarray*}}

\newtheorem{theorem}{Theorem}[section]
\newtheorem{lemma}{Lemma}[section]
\newtheorem{prp}[theorem]{Proposition}
\newtheorem{thm}[theorem]{Theorem}
\newtheorem{cor}[theorem]{Corollary}
\newtheorem{dfn}{Definition}[section]

\begin{document}
\renewcommand{\theequation}{\arabic{section}.\arabic{equation}}
\begin{titlepage}
\title{\bf Large deviation principles  for 3D
stochastic primitive equations
}
\author{Zhao Dong$^{1,}$,\ \ Jianliang Zhai$^{2,}$,\ \ Rangrang Zhang$^{1,}$\\
{\small $^1$ Institute of Applied Mathematics,
Academy of Mathematics and Systems Science, Chinese Academy of Sciences,}\\
{\small No. 55 Zhongguancun East Road, Haidian District, Beijing, 100190, P. R. China}\\
{\small $^2$ School of Mathematical Sciences, University of Science and Technology of China,}\\
 {\small No. 96 Jinzhai Road, Hefei, 230026, P. R. China}\\
({\sf dzhao@amt.ac.cn},\ {\sf zhaijl@ustc.edu.cn},\ {\sf rrzhang@amss.ac.cn} )}
\date{}
\end{titlepage}
\maketitle

\begin{abstract}
 In this paper, we establish the large deviation principle for 3D stochastic primitive equations with small perturbation multiplicative noise. The proof is mainly based on the weak convergence approach.
\end{abstract}

\textbf{Keywords}: Large deviation principle; Laplace principle; Weak convergence approach

\section{Introduction}
The main aim of this paper is to establish large deviation principles (LDP) for 3D stochastic primitive equations, which is a fundamental model in meteorology. In the determined case, the primitive equations are derived from the Navier-Stokes equations, with rotation, coupled with thermodynamics
and salinity diffusion-transport equations, by assuming two important simplifications:
Boussinesq approximation and the hydrostatic balance ( see \cite{L-T-W-1,L-T-W-2,JP}). This model in the determined case has been intensively investigated because of the interests stemmed from physics and mathematics. For example, the mathematical study of the primitive equations originated in a series of articles by J.L. Lions, R. Temam, and S. Wang in the early 1990s \cite{L-T-W-1,L-T-W-2,L-T-W-3,L-T-W-4}, where they set up the mathematical framework and showed the global existence of weak solutions. One remarkable result is that C. Cao and E.S. Titi developed a beautiful approach to dealing with the $L^6$-norm of the fluctuation $\tilde{v}$ of horizontal velocity and obtained the global well-posedness for the 3D viscous primitive equations in  \cite{C-T-1}.

\
For the primitive equations in random case, many results have been obtained. In \cite{Guo}, B. Guo and  D. Huang obtained
the existence of universal random attractor of strong solution under the assumptions that the momentum
equation is driven by an additive stochastic forcing and the
thermodynamical equation is driven by a fixed heat source. A. Debussche, N. Glatt-Holtz, R. Temam and M. Ziane established the global well-posedness of strong solution for the primitive
equations driven by multiplicative random noises in \cite{D-G-T-Z}. In \cite{RR}, the authors obtained the existence of global weak solutions for 3D stochastic primitive equations driven by regular multiplicative noise, and also obtained the exponential mixing property for the weak solutions which are limits of spectral Galerkin approximations. For LDP for stochastic primitive equations, H. Gao and C. Sun obtained a Wentzell-Freidlin type result for the weak solution in \cite{G-S} if this model is driven by small linear multiplicative noise. Moreover, the authors omit the spatial variable $y$ and only take $(x,z)$ into account in order to obtain the global well-posedness of weak solution.
\

In this paper, we consider 3D stochastic primitive equations driven by multiplicative random noise supplied with the same boundary conditions as \cite{D-G-T-Z} and want to establish LDP for its strong solution. As we know, the large deviation theory is concerned with the study of the precise asymptotic behavior governing the decay rate of probabilities of rare events. A classical area
of the large deviation is the Wentzell-Freidlin theory that deals with path probability
asymptotic behavior for stochastic dynamical systems with small noise. 
A weak convergence approach to the
theory of LDP is developed by Dupuis and Ellis in \cite{DE}.  The key idea is to prove some variational representation formula about
the Laplace transform of bounded continuous functionals, which will lead to proving an equivalent
Laplace principle with LDP. In particular, for Brownian functionals, an elegant variational
representation formula has been established by M. Bou\'{e}, P. Dupuis \cite{MP} and A. Budhiraja, Dupuis \cite{BD}.

  \

The proof of small noise LDP is mainly based on the weak convergence approach. Thanks to the equivalence between LDP and the Laplace principle, we only need to verify the Laplace principle holds. A sufficient conditions for the Laplace principle is introduced in Theorem 4.3 of \cite{MP}, which has two parts: the determined part and the random part. During the proof, we focus on the determined part since the random part can be transformed to the determined part. Compared with the primitive equations in \cite{G-S} and 2D geostrophic equations in \cite{Liu}, the difficulty lies in nonlinear terms of our equations is larger since we consider LDP for its strong solution, in that case, $H^1$ estimates is required. Moreover, we can not directly deal with the process that the random solution minus the determined solution and estimate their terms one by one like \cite{Liu} because of the complicate $H^1$ estimates of our equations. Thus, $H^1$ estimates is the key. Fortunately, C. Cao and E.S. Titi developed a beautiful approach to obtain $H^1$ estimates in \cite{C-T-1}, where they consider the fluctuation of horizontal velocity. Based on their work, we obtain the global well-posedness of equation (\ref{equ-9}) by making some additional non-trivial estimates, such as, $|\tilde{v}_h|_{L^{10}(\mathcal{O})}$ estimates and so on. Also, some compact estimates are required. At last, it's worth mentioning that our result is obtained without adding additional regular conditions on the noise, only those in \cite{D-G-T-Z} is enough.

\

This paper is organized as follows. The mathematical formulation for the stochastic primitive equations is in Sects. 2 and 3. Freidlin-Wentzell large deviations and the weak convergence method are introduced in Sect. 4. Then the well-posedness and general a prior estimates for the model are proved in Sect. 5. Finally, a large deviation principle is given in Sect. 6.
\section{Preliminaries}
Let $D$ be a smooth bounded open domain in $\mathbb{R}^2$. Set $\mathcal{O}=D\times (-1,0)$. Consider the
3D primitive equations of the large-scale ocean on $\mathcal{O}\times [0,T]$ driven by a stochastic forcing, in a Cartesian system,
\begin{eqnarray}\label{eq-1}
\frac{\partial v}{\partial t}+(v\cdot \nabla)v+\theta\frac{\partial v}{\partial z}+f{k}\times v +\nabla P +L_1v &=&\psi_1(t,v,T)\frac{dW_1}{dt},\\
\label{eq-2}
\partial_{z}P+T&=&0,\\
\label{eq-3}
\nabla\cdot v+\partial_{z}\theta&=&0, \\
\label{eq-4}
\frac{\partial T}{\partial t}+(v\cdot\nabla)T+\theta\frac{\partial T}{\partial z}+L_2 T&=&\psi_2(t,v,T)\frac{dW_2}{dt},
\end{eqnarray}
where the horizontal velocity field $v=(v_{1},v_{2})$, the three-dimensional velocity field\ $(v_{1},v_{2},\theta)$, the temperature\ $T$ and the pressure\ $P$ are all unknown functionals. $f$ is the Coriolis parameter. ${k}$ is vertical unit vector. $W_1$ and $W_2$ are two independent cylindrical Winner processes which will be given in Sect. \ref{sec 2}. $\nabla=(\partial x,\partial y)$, $\Delta=\partial^{2}_{x}+\partial^{2}_{y}$. The viscosity and the heat diffusion operators $L_1$ and $L_2$ are given by
\begin{eqnarray*}
L_1v&=&-A_h\Delta v -A_v\frac{\partial^2 v}{\partial z^2},\\
L_2T&=&-K_h\Delta T -K_v\frac{\partial^2 T}{\partial z^2},
\end{eqnarray*}
where $A_h$, $A_v$ are positive molecular viscosities and $K_h$, $K_v$ are positive conductivity constants. Without loss of generality, we assume that
$$
A_h=A_v=K_h=K_v=1.
$$

Then, we supply the same boundary conditions as \cite{D-G-T-Z},
\begin{eqnarray}\label{eq1}
 \partial_{z}v=0, \ \theta=0,\  \partial_{z}T=0  & &  \rm {on} \ \textit{D}\times\{0\}=\Gamma_{\textit{u}},\\
 \label{eq2}
 \partial_{z}v=0,\  \theta=0,\  \partial_{z}T=0        & &   \rm {on}\ \textit{D}\times\{-1\}=\Gamma_{\textit{b}},\\
 \label{eq3}
 v=0,\ \frac{\partial T}{\partial {n}}=0  & &   {\rm {on}}\   \partial \textit{D}\times [-1,0]=\Gamma_{l},
 \end{eqnarray}
 where ${n}$ is the normal vector to $\Gamma_{l}$.

Integrating (\ref{eq-3}) from $-1$ to $z$ and using  (\ref{eq1}), (\ref{eq2}), we have
\begin{equation}
\theta(t,x,y,z):=\Phi(v )(t,x,y,z)=-\int^{z}_{-1}\nabla\cdot v (t,x,y,z')dz',
\end{equation}
moreover,
\[
\int^{0}_{-1}\nabla\cdot v  dz=0.
\]
Integrating (\ref{eq-2}) from $-1$ to $z$, set $p_{b}$ be a certain unknown function at $\Gamma_{b}$ satisfying
\[
P(x,y,z,t)= p_{b}(x,y,t)-\int^{z}_{-1} T(x,y,z',t) dz'.
\]
Then, (\ref{eq-1})-(\ref{eq-4}) can be rewritten as
 \begin{eqnarray}\label{eq5-1}
&\frac{\partial v}{\partial t}+(v\cdot \nabla)v+\Phi(v)\frac{\partial v}{\partial z}+f{k}\times v +\nabla p_{b}-\int^{z}_{-1}\nabla T dz' + L_1v =\psi_1(t,v,T)\frac{dW_1}{dt},&\\
\label{eq-6-1}
&\frac{\partial T}{\partial t}+(v\cdot\nabla)T+\Phi(v)\frac{\partial T}{\partial z}+ L_2T=\psi_2(t,v,T)\frac{dW_2}{dt},&\\
\label{eq-7-1}
&\int^{0}_{-1}\nabla\cdot v  dz=0.&
\end{eqnarray}
The boundary value conditions for (\ref{eq5-1})-(\ref{eq-7-1}) are given by
\begin{eqnarray}\label{eq-8-1}
 \partial_{z}v=0,\ \partial_{z}T=0    &&      {\rm on}\ \Gamma_{\textit{u}},\\
 \label{eq-9-1}
 \partial_{z}v=0,\ \partial_{z}T=0     &&   {\rm on}\ \Gamma_{\textit{b}},\\
 \label{eq-10-1}
 v=0,\ \frac{\partial T}{\partial {n}}=0      &&   {\rm on}\ \Gamma_{l}.
 \end{eqnarray}
Denote $Y =(v, T)$ and the initial value conditions are
 \begin{equation}\label{eq-11-1}
 Y(0)=Y_0=(v_{{0}},T_{{0}}).
 \end{equation}
\section{Formulation of this System}\label{sec 2}
\subsection{Some Functional Spaces }
Let $\mathcal{L}(K_1;K_2)$ (resp. $\mathcal{L}_2(K_1;K_2)$) be the space of bounded (resp. Hilbert-Schmidt) linear operators from the Hilbert space $K_1$ to $K_2$, the norm is denoted by $\|\cdot\|_{\mathcal{L}(K_1;K_2)}(\|\cdot\|_{\mathcal{L}_2(K_1;K_2)})$. Denote by $|\cdot|_{L^p(D)}$ the norm of $L^p(D)$ and $|\cdot|_{H^p(D)}$ the norm of $H^p(D)$ for $p\in \mathbb{N}_{+}$. In particular, $|\cdot|$ and $(\cdot,\cdot)$ represent the norm and inner product of $L^2(\mathcal{O})$. For the classical Sobolev space $H^{m}(\mathcal{O})$, $m\in \mathbb{N}_+$,
\begin{equation}\notag
\left\{
  \begin{array}{ll}
    H^{m}(\mathcal{O})=\Big\{U\in (L^2(\mathcal{O}))^3\Big| \partial_{\alpha}U\in L^2(\mathcal{O})\ {\rm for} \ |\alpha|\leq m\Big\},&  \\
    |U|^2_{H^{m}(\mathcal{O})}=\sum_{0\leq|\alpha|\leq m}|\partial_{\alpha}U|^2. &
  \end{array}
\right.
\end{equation}
It's known that $(H^{m}(\mathcal{O}), |\cdot|_{H^{m}(\mathcal{O})})$ is a Hilbert space.

 Define working spaces for the equations (\ref{eq5-1})-(\ref{eq-11-1}). Let
 \begin{eqnarray}\notag
 &&\mathcal{V}_1:=\left\{v\in (C^{\infty}(\mathcal{O}))^2;\  \frac{\partial v}{\partial z}\Big|_{\Gamma_u,\Gamma_b}=0,\  v\Big|_{\Gamma_l}=0,\ \int^{0}_{-1}\nabla \cdot v dz=0\right\},\\ \notag
 &&\mathcal{V}_2:=\left\{T\in C^{\infty}(\mathcal{O});\ \frac{\partial T}{\partial z}\Big|_{\Gamma_u}=0,\ \frac{\partial T}{\partial z}\Big|_{\Gamma_b}=0,\  \frac{\partial T}{\partial {n}}\Big|_{\Gamma_l}=0 \right\},
 \end{eqnarray}
$V_1$= the closure of $ \mathcal{V}_1$ with respect to the norm $|\cdot|_{H^{1}(\mathcal{O})}\times |\cdot|_{H^{1}(\mathcal{O})}$,\\
$V_2$= the closure of $ \mathcal{V}_2$ with respect to the norm $|\cdot|_{H^{1}(\mathcal{O})}$,\\
$H_1$= the closure of $ \mathcal{V}_1$ with respect to the norm $|\cdot|\times|\cdot| $,\\
$H_2$= the closure of $ \mathcal{V}_2$ with respect to the norm $|\cdot|$.

Set \[
V=V_1\times V_2, \quad H=H_1\times H_2.
\]

 The inner products and norms on $V$, $H$ are given by, for any $Y=(v,T), Y_1=(v_1,T_1)$,
\begin{eqnarray*}
&&(Y,Y_1)_{V}=(v,v_1)_{V_1}+(T,T_1)_{V_2},\\
&&(Y,Y_1)=(v,v_1)+(T,T_1)=(v^{(1)},(v_1)^{(1)})+(v^{(2)},(v_1)^{(2)})+(T,T_1),\\
&&\|Y\|_{V}=(Y,Y)^{\frac{1}{2}}_{V}=(v,v)^{\frac{1}{2}}_{V_1}+(T,T)^{\frac{1}{2}}_{V_2}.
\end{eqnarray*}
\subsection{Some Functionals}
Define three bilinear operators $a:V\times V\rightarrow \mathbb{R}$,\ $a_1:V_1\times V_1\rightarrow \mathbb{R}$,\ $a_2:V_2\times V_2\rightarrow \mathbb{R}$,
 and their corresponding linear operators $A: V\rightarrow V^{'}$, $A_1: V_1\rightarrow V^{'}_1$, $A_2: V_2\rightarrow V^{'}_2$ as follows, for any $Y=(v,T)$, $Y_1=(v_1, T_1)\in V$,
 \[
 a(Y,Y_1):=(AY,Y_1)=   
   a_1(v,v_1)+  
  a_2(T,T_1)\\  
  ,
 \]
 where
\begin{eqnarray}\notag
a_1(v,v_1):=(A_1v, v_1)=\int_{\mathcal{O}}\left(\nabla v\cdot \nabla v_1+\frac{\partial v}{\partial z}\cdot\frac{\partial v_1}{\partial z}\right)dxdydz,
\end{eqnarray}
\begin{eqnarray}\notag
a_2(T,T_1):=(A_2T, T_1)=\int_{\mathcal{O}}\left(\nabla T\cdot \nabla T_1+\frac{\partial T}{\partial z}\frac{\partial T_1}{\partial z}\right)dxdydz.
\end{eqnarray}
The following lemma follows Lemma 2.4 in \cite{L-T-W-2}  readily.
\begin{lemma}
\begin{description}
  \item[(i)] The operators $a$, $a_i\ (i=1,2)$ are coercive, continuous, and therefore, the operators $A:V\rightarrow V'$ and $A_i: V_i\rightarrow V'_i\ (i=1,2)$ are isomorphisms. Moreover,
\begin{eqnarray*}
a(Y,Y_1)&\leq& C_1\|Y\|_V\|Y_1\|_V,\\
a(Y,Y)&\geq& C_2\|Y\|^2_V,
\end{eqnarray*}
where $C_1$ and $C_2$ are two positive constants and can be determined in concrete conditions.

  \item[(ii)] The isomorphism  $A:V\rightarrow V'$ (respectively $A_i: V_i\rightarrow V'_i\ (i=1,2)$) can be extended to a self-adjoint unbounded linear  operator on $H$ (respectively on $H_i$, i=1,2), with compact inverse $A^{-1}: H\rightarrow H$ (respectively $A^{-1}_i; H_i\rightarrow H_i \ (i=1,2)$).
\end{description}
\end{lemma}
 It's known that $A_1$ is a self-adjoint operator with discrete spectrum in $H_1$. Denote by $\{k_n\}_{n=1,2,\cdot\cdot\cdot}$ the eigenbasis of  $A_1$ and its associated eigenvalues $\{\nu_n\}_{n=1,2,\cdot\cdot\cdot}$ is increasing. Similarly, $A_2$ is a self-adjoint operator with discrete spectrum in $H_2$. Let $(l_n)_{n=1,2,\cdot\cdot\cdot}$ be the eigenbasis of $A_2$ and its associated increasing eigenvalues $\{\lambda_n\}_{n=1,2,\cdot\cdot\cdot}$. It is easy to see that $\bar{e}_{n,0}=\left(                 
  \begin{array}{c}   
   k_n\\  
   0\\  
  \end{array}
\right)$ and $\bar{e}_{0,m}=\left(                 
  \begin{array}{c}   
   0\\  
   l_m\\  
  \end{array}
\right)$ is the eigenbasis of $(A, {D}(A))$, and we can rearrange $\{\bar{e}_{n,0},\bar{e}_{0,m}\}_{n,m=1,2,\cdots}$, denoted by $\{e_n\}_{n=1,2,\cdots}$, such that the associated eigenvalues is an increasing sequence, denoted by $\{\mu_n\}_{n=1,2,\cdot\cdot\cdot}$.

For any $s\in \mathbb{R}$, the fractional power $(A^s, {D}(A^s))$ of the operator  $(A,{D}(A) )$ is defined as
\begin{equation}\notag
\left\{
  \begin{array}{ll}
    {D}(A^s)=\Big\{Y=\sum_{n=1}^{\infty}y_n e_n \Big|  \sum_{n=1}^{\infty}\mu^{2s}_n|y_n|^2<\infty\Big\}; &  \\
    A^s Y=\sum_{n=1}^{\infty}\mu^s_ny_n e_n, \quad where\  Y=\sum_{n=1}^{\infty}y_n e_n .&
  \end{array}
\right.
\end{equation}
Set
\[
\|Y\|^{A}_s=|A^{\frac{s}{2}}Y|,\quad \mathbb{H}^A_s={D}(A^{\frac{s}{2}}).
\]
It's obvious that $(\mathbb{H}^A_s, \|\cdot\|^A_s)$ is a Hilbert space and $(\mathbb{H}^A_0, \|\cdot\|^A_0)=(H,|\cdot|)$ and $(\mathbb{H}^A_1, \|\cdot\|^A_1)=(V, \|\cdot\|_V)$. For simplicity, denote $\|\cdot\|=\|\cdot\|_V$.
Thanks to the regularity theory of the stokes operator, $\mathbb{H}^A_s$ is a closed subset of $H^s(\mathcal{O})$ and $\|\cdot\|^A_s$ is equivalent to the usual norm $|\cdot|_{H^s(\mathcal{O})}$ for $s\leq 2$.
Similarly, we can define $(\mathbb{H}^{A_1}_s,\|\cdot\|^{A_1}_s)$ and $(\mathbb{H}^{A_2}_s, \|\cdot\|^{A_2}_s)$. For convenience,
all of them will be denoted by $(\mathbb{H}_s,\|\cdot\|_s)$.

Now, we define three mappings $b: V\times V\times V\rightarrow \mathbb{R}$, $b_i: V_1\times V_i\times V_i\rightarrow \mathbb{R}\ (i=1,2)$  and the associated operators $B: V\times V\rightarrow V'$, $B_i: V_1\times V_i\rightarrow V'_i\ (i=1,2)$ by setting
\begin{eqnarray*}
b(Y,Y_1,Y_2)&:=&(B(Y,Y_1),Y_2)=   
   b_1(v,v_1,v_2)+  
   b_2(v, T_1, T_2) 
 ,\\
b_1(v,v_1,v_2)&:=&(B_1(v,v_1), v_2)=\int_{\mathcal{O}}\left[(v\cdot \nabla)v_1+\Phi(v)\frac{\partial v_1}{\partial z}\right]\cdot v_2dxdydz,\\
b_2(v, T_1, T_2)&:=&(B_2(v,T_1), T_2)=\int_{\mathcal{O}}\left[(v\cdot \nabla)T_1+\Phi(v)\frac{\partial T_1}{\partial z}\right]\cdot T_2dxdydz,
\end{eqnarray*}
for any $Y=(v, T)$, $Y_i=(v_i,T_i)\in V$. Then we have

\begin{lemma}\label{lemma-1}
 For any $Y,\ Y_1\in V$,
\[
\left(B(Y,Y_1),Y_1\right)=b(Y,Y_1,Y_1)=b_1(v,v_1,v_1)=b_2(v,T_1,T_1)=0.
\]
\end{lemma}
Moreover, we define another mapping $g: V\times V \rightarrow \mathbb{R}$ and the associated linear operator $G: V\rightarrow V'$ by
\begin{eqnarray}\notag
g(Y,Y_1)&:=&(G(Y), Y_1)\notag
\\      &=&\int_{\mathcal{O}}\left[f(k\times v)\cdot v_1+(\nabla p_b-\int^z_{-1}\nabla Tdz')\cdot v_1\right]dxdydz 
. \notag
\end{eqnarray}
By (\ref{eq-7-1}), we have
\begin{equation}\notag
(v,\nabla p_{b})=\left(\int^{0}_{-1}vdz,\nabla p_b\right)_{L^2(D)}=-\left(p_b,\int^{0}_{-1}\nabla \cdot vdz \right)_{L^2(D)}=0,
\end{equation} and by $(v,f{k}\times v)=0$,
we have
\begin{lemma}\label{lemma-2}
\begin{description}\notag
  \item[(i)]
\[
g(Y,Y)=(G(Y),Y)=
-\int_{\mathcal{O}}\Big[\Big(\int^z_{-1}\nabla Tdz'\Big)\cdot v\Big]dxdydz
.
\]
  \item[(ii)] There exists a constant $C$, such that
\begin{eqnarray}\label{a-2}
|(G(Y),Y)|&\leq& C(|T|\|v\|_V\vee\|T\|_V|v|),\\
\label{a-3}
|(G(Y),Y_1)|&\leq& C|v||v_1|+C(|T|\|v_1\|_V\vee\|T\|_V|v_1|).
\end{eqnarray}
\end{description}
\end{lemma}

Using the functionals defined above, we merge (\ref{eq5-1}) and (\ref{eq-6-1}) as follows
\begin{eqnarray}\label{equ-7}
\left\{
  \begin{array}{ll}
    dY(t)+AY(t)dt+B(Y(t),Y(t))dt+G(Y(t))dt=\psi(t,Y(t)) dW(t), \\
    Y(0)=Y_0.
  \end{array}
\right.
\end{eqnarray}
where
\begin{equation}\notag
W=\left(                 
  \begin{array}{c}   
    W_1\\  
    W_2 \\  
  \end{array}
\right) ,\quad
\psi(t,Y(t))
=\left(                 
  \begin{array}{cc}   
   \psi_1(t,Y(t)) & 0 \\  
   0 & \psi_2(t,Y(t)) \\  
  \end{array}
\right).
\end{equation}

\subsection{Some Inequalities}
Let us recall some interpolation inequalities used later (see Sect. 4.1 in \cite{Guo}).\\
For $h\in H^1(D)$,
\begin{eqnarray*}
|h|_{L^4(D)}&\leq& c|h|^{\frac{1}{2}}_{L^2(D)}|h|^{\frac{1}{2}}_{H^1(D)},\\
|h|_{L^5(D)}&\leq& c|h|^{\frac{3}{5}}_{L^3(D)}|h|^{\frac{2}{5}}_{H^1(D)},\\
|h|_{L^6(D)}&\leq& c|h|^{\frac{2}{3}}_{L^4(D)}|h|^{\frac{1}{3}}_{H^1(D)}.
\end{eqnarray*}
For $h\in H^1(\mathcal{O})$,
\begin{eqnarray*}
|h|_3&\leq& c|h|^{\frac{1}{2}}|h|^{\frac{1}{2}}_{H^1(\mathcal{O})},\\
|h|_4&\leq& c|h|^{\frac{1}{4}}|h|^{\frac{3}{4}}_{H^1(\mathcal{O})},\\
|h|_6&\leq& c|h|_{H^1(\mathcal{O})},\\
|h|_{\infty}&\leq& c|h|^{\frac{1}{2}}_{H^1(\mathcal{O})}|h|^{\frac{1}{2}}_{H^2(\mathcal{O})}.
\end{eqnarray*}
Using the similar argument as page 17 in \cite{C-T-1} and Proposition 2.2 in \cite{C-T-2}, we have
\begin{lemma}\label{le-2}
Let $u, f,g $ be smooth functions, then
\begin{description}
 \item[(i)] $|\int_{\mathcal{O}}g\cdot[(u \cdot \nabla)f]dxdydz|\leq c|\nabla f||g|_3|u|_6\leq c |\nabla f||g|^{\frac{1}{2}}|\nabla g|^{\frac{1}{2}}|\nabla u|$,
 \item[(ii)] $|\int_{\mathcal{O}}\Phi(u)f\cdot g dxdydz|\leq c|\nabla u||g|^{\frac{1}{2}}|\nabla g|^{\frac{1}{2}}|f|^{\frac{1}{2}}|\nabla f|^{\frac{1}{2}}$,
     \item[(iii)] $|\int_{\mathcal{O}}\Phi(u)f\cdot g dxdydz|\leq c|f||\nabla u|^{\frac{1}{2}}\| u\|^{\frac{1}{2}}_2|\nabla g|^{\frac{1}{2}}| g|^{\frac{1}{2}}$.
\end{description}
\end{lemma}
At last, we recall the integral version of Minkowshy inequality for the $L^p$ spaces, $p\geq 1$. Let $\mathcal{O}_1\subset \mathbb{R}^{m_1}$ and $\mathcal{O}_2\subset \mathbb{R}^{m_2}$ be two Borel measurable subsets, where $m_1$ and $m_2$ are two positive integers. Suppose that $f(\xi, \eta)$ is measurable over $\mathcal{O}_1\times \mathcal{O}_2$. Then
\begin{equation}\notag
\left[\int_{\mathcal{O}_1}\left(\int_{\mathcal{O}_2}|f(\xi,\eta)|d\eta\right)^p d\xi\right]^{1/p}\leq \int_{\mathcal{O}_2}\left(\int_{\mathcal{O}_1}|f(\xi, \eta)|^p d\xi\right)^{1/p}d\eta.
\end{equation}

\subsection{Definition of Strong Solution}
For the strong solution of (\ref{equ-7}), we shall fix a single stochastic basis $\mathcal{T}:=(\Omega, \mathcal{F}, \{\mathcal{F}_t\}_{t\geq 0}, \mathbb{P}, W)$. Here,
\[
W=\left(                 
  \begin{array}{c}   
    W_1\\  
    W_2 \\  
  \end{array}
\right)
 \]
 is a cylindrical Brownian motion with the form $W(t,\omega)=\sum_{i\geq1}r_iw_i(t,\omega)$, where $\{r_i\}_{i\geq 1}$ is a complete orthonormal basis of a Hilbert space
 \[
U=\left(                 
  \begin{array}{c}   
    U_1\\  
    U_2 \\  
  \end{array}
\right)
 \]
 and $\{w_i\}_{i\geq1}$ is a sequence of independent one-dimensional standard Brownian motions on $(\Omega, \mathcal{F}, \{\mathcal{F}_t\}_{t\geq 0}, \mathbb{P})$, $U_1$ and $U_2$ are separable Hilbert spaces.

Given any pair of Banach spaces $\mathcal{X}$ and $\mathcal{Y}$,  $Bnd_u(\mathcal{X}, \mathcal{Y})$ stands for the collection of all continuous mappings $\psi: [0,\infty)\times \mathcal{X} \rightarrow \mathcal{Y}$ such that
\[
\|\psi(t,x)\|_{\mathcal{Y}}\leq c(1+\|x\|_{\mathcal{X}}), \quad x\in \mathcal{X}, \ t\geq 0,
\]
where the numerical constant $c$ may be chosen independent of $t$. If,  in addition,
\[
\|\psi(t,x)-\psi(t,y)\|_{\mathcal{Y}}\leq c\|x-y\|_{\mathcal{X}}, \quad x,y\in \mathcal{X}, \ t\geq 0,
\]
we say $\psi$ is in $Lip_u(\mathcal{X},\mathcal{Y})$.
\begin{description}
  \item[\textbf{Hypothesis H0}] We assume that $\psi: [0,\infty)\times H\rightarrow \mathcal{L}_2({U}, H)$ with
\begin{eqnarray*}
\psi\in Lip_u(H, \mathcal{L}_2({U}; H))\cap Lip_u(V, \mathcal{L}_2({U}; V))\cap Bnd_u(V, \mathcal{L}_2({U}; D(A))).
\end{eqnarray*}
\end{description}
Now, we give the definition of strong solution to (\ref{equ-7}).
\begin{dfn}\label{dfn-3}\cite{D-G-T-Z}
Let $\mathcal{T}=(\Omega, \mathcal{F}, \{\mathcal{F}_t\}_{t\geq 0}, \mathbb{P}, W)$ be a fixed stochastic basis and suppose that $Y_0\in V$.
$Y$ is called a strong solution of (\ref{equ-7}) if $Y(\cdot)$ is an $\mathcal{F}_t-$ adapted process in $V$,  such that
\[
Y(\cdot)\in L^2(\Omega; C([0,T];V))\bigcap L^2(\Omega; L^2([0,T];D(A))),\ \ \forall T>0,
\]
and for every $t\geq 0$,
\[
 Y(t)+\int^{t}_{0}\Big(AY+B(Y,Y)+G(Y)\Big)ds=Y_0+\int^{t}_{0}\psi(s,Y(s)) dW(s),
 \]
 holds in $V'$, $\mathbb{P}-$ a.s.
\end{dfn}
\begin{thm}\label{thm-3}\cite{D-G-T-Z}
Suppose that $Y_0\in V$. Assume that \textbf{Hypothesis H0} holds for $\psi$. Then there exists a unique global solution $Y$ of (\ref{equ-7}) in the sense of Definition  \ref{dfn-3} with $Y(0)=Y_0$.
\end{thm}
\section{Freidlin-Wentzell's Large Deviations}
In this section, we consider the large deviation principle for the stochastic primitive. Here, we will use the weak convergence approach introduced by Budhiraja and Dupuis in \cite{BD}. Let us first recall some standard definitions and results from large deviation theory (see \cite{DZ})

Let $\{Y^\varepsilon\}$ be a family random variables defined on a probability space $(\Omega, \mathcal{F}, \mathbb{P})$ taking values in some Polish space $\mathcal{E}$.

\begin{dfn}
(Rate Function) A function $I: \mathcal{E}\rightarrow [0,\infty]$ is called a rate function if $I$ is lower semicontinuous. A rate function $I$ is called a good rate function if the level set $\{x\in \mathcal{E}: I(x)\leq M\}$ is compact for each $M<\infty$.
\end{dfn}
\begin{dfn}
\begin{description}
  \item[(i)] (Large deviation principle) The sequence $\{Y^\varepsilon\}$ is said to satisfy the large deviation principle with rate function $I$ if for each Borel subset $A$ of $\mathcal{E}$
      \[
      -\inf_{x\in A^o}I(x)\leq \lim \inf_{\varepsilon\rightarrow 0}\varepsilon \log \mathbb{P}(Y^\varepsilon\in A)\leq \lim \sup_{\varepsilon\rightarrow 0}\varepsilon \log \mathbb{P}(Y^\varepsilon\in A)\leq -\inf_{x\in \bar{A}}I(x),
      \]
      where $A^o$ and $\bar{A}$ denote the interior and closure of $A$ in $\mathcal{E}$, respectively.
  \item[(ii)] (Laplace principle) The sequence $\{Y^\varepsilon\}$ is said to satisfy the Laplace principle with rate function $I$ if for each bounded continuous real-valued function $f$ defined on $\mathcal{E}$
      \[
      \lim_{\varepsilon\rightarrow 0}\varepsilon \log E\Big\{\exp[-\frac{1}{\varepsilon}f(Y^\varepsilon)]\Big\}=-\inf_{x\in \mathcal{E}}\{f(x)+I(x)\}.
      \]
\end{description}
\end{dfn}
It well-known that the large deviation principle and the Laplace principle are equivalent if $\mathcal{E}$ is a Polish space and the rate function is good. The equivalence is essentially a consequence of Varadhan's
lemma and Bryc's converse theorem (see \cite{DZ}).

Suppose $W(t)$ is a cylindrical Wiener process on a Hilbert space $U$ defined on a probability space $(\Omega, \mathcal{F},\{\mathcal{F}_t\}_{t\in [0,T]}, \mathbb{P} )$ ( the paths of $W$ take values in $C([0,T];\mathcal{U})$, where $\mathcal{U}$ is another Hilbert space such that the embedding $U\subset \mathcal{U}$ is Hilbert-Schmidt).
Now we define
\begin{eqnarray*}
&\mathcal{A}=\{\phi: \phi\ is\ a\ U\text{-}valued\ \{\mathcal{F}_t\}\text{-}predictable\ process\ s.t.\ \int^T_0 |\phi(s)|^2_Uds<\infty\ a.s.\};\\
&T_M=\{ h\in L^2([0,T];U): \int^T_0 |h(s)|^2_Uds\leq M\};\\
&\mathcal{A}_M=\{\phi\in \mathcal{A}: \phi(\omega)\in T_M,\ \mathbb{P}\text{-}a.s.\}.
\end{eqnarray*}
Here, we use the weak topology on the set $T_M$ under which $T_M$ is a compact space.

Suppose $\mathcal{G}^{\varepsilon}: C([0,T];U)\rightarrow \mathcal{E}$ is a measurable map and $Y^{\varepsilon}=\mathcal{G}^{\varepsilon}(W)$. Now, we list the following sufficient conditions for the Laplace principle (equivalently, large deviation principle) of $Y^{\varepsilon}$ as $\varepsilon\rightarrow 0$.
\begin{description}
  \item[\textbf{Hypothesis H1} ] There exists a measurable map $\mathcal{G}^0: C([0,T];U)\rightarrow \mathcal{E}$ such that the following conditions hold
\end{description}
\begin{description}
  \item[(i)] For every $M<\infty$, let $\{h_{\varepsilon}: \varepsilon>0\}$ $\subset \mathcal{A}_M$. If $h_{\varepsilon}$ converges to $h$ as $T_M$-valued random elements in distribution, then $\mathcal{G}^{\varepsilon}(W(\cdot)+\frac{1}{\sqrt{\varepsilon}}\int^{\cdot}_{0}h_\varepsilon(s)ds)$ converges in distribution to $\mathcal{G}^0(\int^{\cdot}_{0}h(s)ds)$.
  \item[(ii)] For every $M<\infty$, the set $K_M=\{\mathcal{G}^0(\int^{\cdot}_{0}h(s)ds): h\in T_M\}$ is compact subset of $\mathcal{E}$.
\end{description}
\begin{thm}\label{thm-2}
If $\{\mathcal{G}^{\varepsilon}\}$ satisfies \textbf{Hypothesis H1}, then $Y^{\varepsilon}$ satisfies the Laplace principle (hence large deviation principle) on $\mathcal{E}$ with the good rate function $I$ given by
\begin{eqnarray}\label{eq-5}
I(f)=\inf_{\{h\in L^2([0,T];U): f= \mathcal{G}^0(\int^{\cdot}_{0}h(s)ds)\}}\Big\{\frac{1}{2}\int^T_0|h(s)|^2_{U}ds\Big\},\ \ \forall f\in\mathcal{E}.
\end{eqnarray}
By convention, $I(f)=\infty$, if  $\Big\{h\in L^2([0,T];U): f= \mathcal{G}^0(\int^{\cdot}_{0}h(s)ds)\Big\}=\emptyset.$
\end{thm}
\section{Prior Estimates}\label{Sec 4}
Consider the 3D stochastic primitive equations driven by small multiplicative noise
\begin{eqnarray}\label{equ-8}
\left\{
  \begin{array}{ll}
    dY^\varepsilon(t)+AY^\varepsilon(t)dt+B(Y^\varepsilon(t),Y^\varepsilon(t))dt+G(Y^\varepsilon(t))dt
    =\sqrt{\varepsilon}\psi(t,Y^\varepsilon)dW(t), \\
    Y^\varepsilon(0)=Y_0,
  \end{array}
\right.
\end{eqnarray}
where $Y_0\in V$. Under \textbf{Hypothesis H0}, by Theorem \ref{thm-3}, there exists a pathwise unique strong solution of (\ref{equ-8}) in $\Re:= C([0,T];V)\cap L^2([0,T]; D(A)) $, the norm in $\Re$ is that
\[
|Y|^2_{\Re}:=\sup_{0\leq t\leq T}\|Y(t)\|^2+\int^T_0\|Y(t)\|^2_{D(A)}dt.
\]
Therefore, there exist Borel-measurable functions
\begin{eqnarray}\label{eq def G epsilon}
\mathcal{G}^{\varepsilon}: C([0,T];U)\rightarrow \Re\text{  such that }Y^{\varepsilon}(\cdot)=\mathcal{G}^{\varepsilon}(W(\cdot)).
\end{eqnarray}

Now, the aim is to prove the large deviation principle for $Y^{\varepsilon}$.

For $h\in L^2([0,T];U)$, we consider the following skeleton equation
\begin{eqnarray}\label{equ-9}
\left\{
  \begin{array}{ll}
    dY_h(t)+AY_h(t)dt+B(Y_h(t),Y_h(t))dt+G(Y_h(t))dt=\psi(t,Y_h(t)) h(t)dt, \\
    Y_h(0)=Y_0.
  \end{array}
\right.
\end{eqnarray}
Denote by $h=(h_1,h_2)$, we rewrite (\ref{equ-9}) as
\begin{eqnarray} \label{eq-47}
dv_h+[(v_h\cdot \nabla)v_h+\Phi(v_h)\frac{\partial v_h}{\partial z}]dt+(f{k}\times v_h +\nabla p_{b}-\int^{z}_{-1}\nabla T_h dz')dt + L_1v_h dt =\psi_1(t,Y_h)h_1(t)dt,\\
\label{eq-48}
dT_h+[(v_h\cdot\nabla)T_h+\Phi(v_h)\frac{\partial T_h}{\partial z}]dt+ L_2 T_hdt=\psi_2(t,Y_h)h_2(t)dt.
\end{eqnarray}

\subsection{Global Well-posedness }
\begin{thm}\label{thm-4}
Assume \textbf{Hypothesis H0} holds and the initial data $Y_0=(v_0,T_0)\in V$, let $h\in T_M $, then for any $T>0$, (\ref{equ-9}) has a unique strong solution $Y_h\in C([0,T];V)\bigcap L^2([0,T];D(A))$ on the interval $[0,T]$, which depends continuously on the initial data.
\end{thm}
In order to prove Theorem \ref{thm-4}, we need to repeat and partial refined some calculations in  \cite{C-T-1}.

\subsubsection{A priori estimates in H}
Taking the inner product of the equation (\ref{equ-9}) with  $Y_{h}$ in $L^2(\mathcal{O})$, we get
\begin{eqnarray*}
\frac{1}{2}d|Y_{h}|^2+(|\nabla Y_{h}|^2+|\partial_z Y_{h}|^2)dt=-(B(Y_h,Y_h),Y_h)dt-(G(Y_h),Y_h)dt+(\psi(t,Y_h) h,Y_h)dt,
\end{eqnarray*}
by Lemma \ref{lemma-1} and Lemma \ref{lemma-2},
\[
\frac{1}{2}d|Y_{h}|^2+(|\nabla Y_{h}|^2+|\partial_z Y_{h}|^2)dt
\leq
C|Y_{h}|\|Y_h\|dt+C|Y_h||\psi(t,Y_h) h|dt,
\]
by H\"{o}lder inequality and the Young inequality, we have
\[
\frac{1}{2}d|Y_{h}|^2+(|\nabla Y_{h}|^2+|\partial_z Y_{h}|^2)dt\leq
\varepsilon \|Y_{h}\|^2dt+C|Y_h|^2dt+C|\psi(t,Y_h) h|^2dt.
\]
It follows from \textbf{Hypothesis H0} that
\begin{eqnarray}\label{equ-26}
|\psi(t,Y_h) h|^2 &\leq& \|\psi(t,Y_h)\|^2_{\mathcal{L}_2(U;H)}|h|^2_{U}\\ \notag
&\leq& C(1+|Y_h|^2)|h|^2_{U},
\end{eqnarray}
then,
\[
d|Y_{h}|^2+\|Y_{h}\|^2dt\leq C(1+|h|^2_{U})|Y_h|^2 dt +C|h|^2_{U}dt.
\]
Applying Gronwall inequality, we have
\begin{eqnarray}
\sup_{t\in[0,T]}|Y_{h}(t)|^2\leq C_1(|Y_0|^2,M),
\end{eqnarray}
and
\begin{eqnarray}\label{equ-12}
\sup_{t\in[0,T]}|Y_{h}(t)|^2+\int^T_0\|Y_{h}(t)\|^2dt\leq K_1(|Y_0|^2,M),
\end{eqnarray}
where
\begin{eqnarray*}
C_1(|Y_0|^2,M)&=&C(1+M)e^{C(1+M)}(|Y_0|^2+CM),\\
K_1(|Y_0|^2,M)&=&C(1+M)^2e^{C(1+M)}(|Y_0|^2+CM).
\end{eqnarray*}
\subsubsection{Splitting}
From now on, keeping in mind that we consider the case $\alpha=0$ and the model is supplied with and boundary conditions (\ref{eq-8-1})--(\ref{eq-10-1}) in \cite{C-T-1}, let
\[
\bar{v}_h(x,y,t)=\int^0_{-1}v_h(x,y,z',t)dz', \ and\ the\ fluctuation\ \tilde{v}_h=v_h-\bar{v}_h, \quad h=(h_1,h_2),
\]
refer to equation (32) in \cite{C-T-1}, we obtain
\begin{eqnarray}\label{equ-10}
&\frac{\partial \bar{v}_h}{\partial t}-\Delta \bar{v}_h+(\bar{v}_h\cdot \nabla)\bar{v}_h+\overline{[(\tilde{v}_h\cdot\nabla)\tilde{v}_h+(\nabla\cdot\tilde{v}_h)\tilde{v}_h]}
+\nabla p_s(x,y,t)+fk\times \bar{v}_h\\ \notag
\quad\quad &-\nabla[\int^0_{-1}\int^z_{-1}T_h(x,y,z',t)dz'dz]=\int^0_{-1}\psi_1(t,Y_h(t))h_1(t)dz,\\ \notag
&\nabla\cdot \bar{v}_h=0, \ in \ D,\\ \notag
&\bar{v}_h=0,\ on\  \partial D.\notag
\end{eqnarray}
By subtracting (\ref{equ-10}) from (\ref{eq-47}),  $\tilde{v}_h$  satisfies
\begin{eqnarray}\label{equ-11}
&\frac{\partial \tilde{v}_h}{\partial t}+L_1\tilde{v}_h+(\tilde{v}_h\cdot \nabla)\tilde{v}_h-(\int^z_{-1}\nabla\cdot \tilde{v}_h(x,y,z',t)dz')\frac{\partial \tilde{v}}{\partial z}
+(\tilde{v}_h\cdot \nabla)\bar{v}_h+(\bar{v}_h\cdot\nabla)\tilde{v}_h+fk\times \tilde{v}_h\\ \notag
&-\overline{[(\tilde{v}_h\cdot \nabla)\tilde{v}_h+(\nabla\cdot \tilde{v}_h)\tilde{v}_h]}-\nabla\Big(\int^z_{-1}T_h(x,y,z',t)dz'-\int^0_{-1}\int^z_{-1}T_h(x,y,z',t)dz'dz\Big)\\ \notag
&=\psi_1(t,Y_h(t))h_1(t)-\int^0_{-1}\psi_1(t,Y_h(t))h_1(t)dz,\\ \notag
&\frac{\partial \tilde{v}_h}{\partial z}|_{z=0}=0, \ \frac{\partial \tilde{v}_h}{\partial z}|_{z=-1}=0,\ \tilde{v}_h\cdot n|_{\Gamma_l}=0,\ \tilde{v}_h|_{\Gamma_l}=0.
\end{eqnarray}
\subsubsection{$H^1$ estimates}
\textbf{$L^6$ estimate of $\tilde{v}_h$. } \quad Taking the inner product of (\ref{equ-11}) with $|\tilde{v}_h|^4\tilde{v}_h$ in $L^2(\mathcal{O})$. In the same way as Page 10 in \cite{C-T-1}, we obtain
\begin{eqnarray}\label{eq P10 eq1}
&&\frac{d|\tilde{v}_h|^6}{dt}+2\int_{\mathcal{O}}(|\tilde{v}_h|^2|\nabla|\tilde{v}_h|^2|^2+|\tilde{v}_h|^4|\nabla \tilde{v}_h|^2)dxdydz+2\int_{\mathcal{O}}(|\tilde{v}_h|^2|\partial_z|\tilde{v}_h|^2|^2+|\tilde{v}_h|^4|\partial_z \tilde{v}_h|^2)dxdydz\nonumber\\
&\leq& C|\bar{v}_h|^2|\nabla\bar{v}_h|^2|\tilde{v}_h|^6+C|\tilde{v}_h|^6|\nabla\tilde{v}_h|^2+C|\bar{T}_h|^2|\nabla \bar{T}_h|^2+C|\tilde{v}_h|^2|\tilde{v}_h|^6\nonumber\\
&& \quad +\Big|\int_{\mathcal{O}}\Big(\psi_1h_1(t)-\int^0_{-1}\psi_1h_1(t)dz\Big)\cdot |\tilde{v}_h|^4 \tilde{v}_h dxdydz\Big|,
\end{eqnarray}
we only need to estimate the following additional term,
\begin{eqnarray*}
&& \Big|\int_{\mathcal{O}}\Big(\psi_1h_1(t)-\int^{0}_{-1}\psi_1h_1(t)dz\Big)\cdot |\tilde{v}_h|^4 \tilde{v}_h dxdydz\Big|\\
&\leq&\Big( \int_{\mathcal{O}}|\psi_1h_1(t)-\int^0_{-1}\psi_1h_1(t)dz|^2dxdydz\Big)^{\frac{1}{2}}\Big( \int_{\mathcal{O}}|\tilde{v}_h|^{10}dxdydz\Big)^{\frac{1}{2}}\\
&:=& I_1I_2,
\end{eqnarray*}
where H\"{o}lder inequality is used.

For the first term $I_1$, by \textbf{Hypothesis H0} and (\ref{equ-12}), we have
\begin{eqnarray*}
 I_1
&\leq& C|\psi_1(t,Y_h)h_1(t)|\\
&\leq & C(1+|Y_h(t)|)|h_1(t)|_{U}\\
&\leq & C(1+\sup_{t\in[0,T]}|Y_h(t)|)|h_1(t)|_{U}\\
&\leq & C|h_1(t)|_{U}.
\end{eqnarray*}
For the second term $I_2$, by Sobolev inequality, $|u|_{L^{\frac{10}{3}}(\mathcal{O})}\leq C\|u\|^{\frac{3}{5}}|u|^{\frac{2}{5}}_{L^2(\mathcal{O})}$, we have
\begin{eqnarray*}
|\tilde{v}_h|^{10}_{L^{10}(\mathcal{O})}&=&||\tilde{v}_h|^3|^{\frac{10}{3}}_{L^{\frac{10}{3}}(\mathcal{O})}\\
&\leq& C\||\tilde{v}_h|^3\|^2||\tilde{v}_h|^3|^{\frac{4}{3}}_{L^2(\mathcal{O})}\\
&\leq& C||\tilde{v}_h|^3|^{\frac{4}{3}}_{L^2(\mathcal{O})}(||\tilde{v}_h|^3|^2+|\nabla|\tilde{v}_h|^3|^2+|\partial_z|\tilde{v}_h|^3|^2)\\
&\leq& C|\tilde{v}_h|^4_{L^6(\mathcal{O})}\Big[|\tilde{v}_h|^6_{L^6(\mathcal{O})}+\int_{\mathcal{O}}(|\tilde{v}_h|^2|\nabla|\tilde{v}_h|^2|^2+|\tilde{v}_h|^4|\nabla \tilde{v}_h|^2)dxdydz+\int_{\mathcal{O}}(|\tilde{v}_h|^2|\partial_z|\tilde{v}_h|^2|^2+|\tilde{v}_h|^4|\partial_z \tilde{v}_h|^2)dxdydz\Big]\\
&\leq& C|\tilde{v}_h|^{10}_{L^6(\mathcal{O})}+C|\tilde{v}_h|^4_{L^6(\mathcal{O})}\Big[\int_{\mathcal{O}}(|\tilde{v}_h|^2|\nabla|\tilde{v}_h|^2|^2+|\tilde{v}_h|^4|\nabla \tilde{v}_h|^2)dxdydz+\int_{\mathcal{O}}(|\tilde{v}_h|^2|\partial_z|\tilde{v}_h|^2|^2+|\tilde{v}_h|^4|\partial_z \tilde{v}_h|^2)dxdydz\Big],
\end{eqnarray*}
then we obtain
\begin{eqnarray*}
&& \Big|\int_{\mathcal{O}}\Big(\psi_1h_1(t)-\int^0_{-1}\psi_1h_1(t)dz\Big)\cdot |\tilde{v}_h|^4 \tilde{v}_h dxdydz\Big|\\
&\leq&  C|h_1(t)|_{U}|\tilde{v}_h|^{5}_{L^6(\mathcal{O})}+C|h_1(t)|_{U}|\tilde{v}_h|^2_{L^6(\mathcal{O})}
\Big[\int_{\mathcal{O}}(|\tilde{v}_h|^2|\nabla|\tilde{v}_h|^2|^2+|\tilde{v}_h|^4|\nabla \tilde{v}_h|^2)dxdydz\Big]^{\frac{1}{2}}
\\
&&\quad +C|h_1(t)|_{U}|\tilde{v}_h|^2_{L^6({\mathcal{O}})}\Big[\int_{\mathcal{O}}(|\tilde{v}_h|^2|\partial_z|\tilde{v}_h|^2|^2+|\tilde{v}_h|^4|\partial_z \tilde{v}_h|^2)dxdydz\Big]^{\frac{1}{2}}\\
&:=& I_3+I_4+I_5,
\end{eqnarray*}
for $I_3$,
\begin{eqnarray*}
I_3
&\leq & C(1+|\tilde{v}_h|^6_{L^6(\mathcal{O})})(1+|h_1(t)|^2_{U})\\
&\leq & C|h_1(t)|^2_{U}|\tilde{v}_h|^6_{L^6(\mathcal{O})}+C|\tilde{v}_h|^6_{L^6(\mathcal{O})}+C|h_1(t)|^2_{U}+C,
\end{eqnarray*}
for $I_4$,
\begin{eqnarray*}
I_4
&\leq &
\varepsilon\int_{\mathcal{O}}(|\tilde{v}_h|^2|\nabla|\tilde{v}_h|^2|^2+|\tilde{v}_h|^4|\nabla \tilde{v}_h|^2)dxdydz+C|h_1(t)|^2_{U}|\tilde{v}_h|^4_{L^6(\mathcal{O})}\\
&\leq & \varepsilon\int_{\mathcal{O}}(|\tilde{v}_h|^2|\nabla|\tilde{v}_h|^2|^2+|\tilde{v}_h|^4|\nabla \tilde{v}_h|^2)dxdydz+C|h_1(t)|^2_{U}|\tilde{v}_h|^6_{L^6(\mathcal{O})}+C|h_1(t)|^2_{U},
\end{eqnarray*}
$I_5$ is similar to $I_4$,
\begin{eqnarray*}
I_5
&\leq &
\varepsilon\int_{\mathcal{O}}(|\tilde{v}_h|^2|\partial_z|\tilde{v}_h|^2|^2+|\tilde{v}_h|^4|\partial_z \tilde{v}_h|^2)dxdydz+C|h_1(t)|^2_{U}|\tilde{v}_h|^4_{L^6(\mathcal{O})}\\
&\leq & \varepsilon\int_{\mathcal{O}}(|\tilde{v}_h|^2|\partial_z|\tilde{v}_h|^2|^2+|\tilde{v}_h|^4|\partial_z \tilde{v}_h|^2)dxdydz+C|h_1(t)|^2_{U}|\tilde{v}_h|^6_{L^6(\mathcal{O})}+C|h_1(t)|^2_{U},
\end{eqnarray*}
thus, we obtain
\begin{eqnarray}\label{equ-27}
&& |\int_{\mathcal{O}}(\psi_1h_1(t)-\int^0_{-1}\psi_1h_1(t)dz)\cdot |\tilde{v}_h|^4 \tilde{v}_h dxdydz|\\ \notag
&\leq& \varepsilon\int_{\mathcal{O}}(|\tilde{v}_h|^2|\nabla|\tilde{v}_h|^2|^2+|\tilde{v}_h|^4|\nabla \tilde{v}_h|^2)dxdydz+\varepsilon\int_{\mathcal{O}}(|\tilde{v}_h|^2|\partial_z|\tilde{v}_h|^2|^2+|\tilde{v}_h|^4|\partial_z \tilde{v}_h|^2)dxdydz\\ \notag
&&\quad +C(1+|h_1(t)|^2_{U})|\tilde{v}_h|^6_{L^6(\mathcal{O})}
+C(1+|h_1(t)|^2_{U}).
\end{eqnarray}
Putting (\ref{equ-12}), (\ref{eq P10 eq1}) and (\ref{equ-27}) together, we have
\begin{eqnarray}\notag
|\tilde{v}_h(t)|^6_{L^6(\mathcal{O})}&+&\int^t_0\Big(\int_{\mathcal{O}}(|\tilde{v}_h|^2|\nabla|\tilde{v}_h|^2|^2+|\tilde{v}_h|^4|\nabla \tilde{v}_h|^2)dxdydz\\ &+&\int_{\mathcal{O}}(|\tilde{v}_h|^2|\partial_z|\tilde{v}_h|^2|^2+|\tilde{v}_h|^4|\partial_z \tilde{v}_h|^2)dxdydz\Big)ds
\leq K_2(t),\label{equ-13}
\end{eqnarray}
where
\[
K_2(t)=e^{\Big(C(1+M)K^2_1(t)\Big)}\Big[\|v_0\|^6+C(1+M)+K^2_1(t)\Big].
\]
\textbf{$L^6$ estimates of $T_h$.} \quad  It's similar to $L^6$ estimates of $\tilde{v}_h$, we obtain
\begin{eqnarray}\label{eq-50}
|T_h(t)|^6_{L^6(\mathcal{O})}+\int^t_0\Big(\int_{\mathcal{O}}|T_h|^4|\nabla T_h|^2dxdydz+\int_{\mathcal{O}}|T_h|^4|\frac{\partial T_h}{\partial z}|^2dxdydz\Big)ds
\leq K_3(t),
\end{eqnarray}
where
\[
K_3(t)=e^{\Big(C(1+M)\Big)}\Big[\|T_0\|^6+C(1+M)\Big].
\]
\textbf{$|\nabla \bar{v}_h|_{L^2(D)}$ estimates.} \quad Taking the inner product of equation (\ref{equ-10}) with $-\Delta \bar{v}_h $ in $L^2(D)$,  as Page 12  in \cite{C-T-1},
we have
\begin{eqnarray}\label{eq P12 star}
&&\frac{d|\nabla \bar{v}_h|^2}{dt}+2|\Delta \bar{v}_h|^2\\
&\leq& C|\bar{v}_h|^2|\nabla \bar{v}_h|^4+C|\nabla\tilde{v}_h|^2+C\int_{\mathcal{O}}|\tilde{v}_h|^4|\nabla \tilde{v}_h|^2dxdydz+C|\bar{v}_h|^2+|\int_{D}\Delta\bar{v}_h\Big(\int^0_{-1}\psi_1(t,Y_h)h_1(t)dz'\Big)dxdy|,\nonumber
\end{eqnarray}
 we only need to estimate the following additional term
\begin{eqnarray*}
|\int_{D}\Delta\bar{v}_h\Big(\int^0_{-1}\psi_1(t,Y_h)h_1(t)dz'\Big)dxdy|
&\leq& C|\Delta\bar{v}_h|_{L^2(D)}|\int^0_{-1}\psi_1(t,Y_h)h_1(t)dz'|_{L^2(D)}\\
&\leq& \varepsilon |\Delta\bar{v}_h|^2_{L^2(D)}+C|\int^0_{-1}\psi_1(t,Y_h)h_1(t)dz'|^2_{L^2(D)}.
\end{eqnarray*}
Since, by (\ref{equ-12}) and $\rm H\ddot{o}lder's$ inequality,
\begin{eqnarray*}
|\int^0_{-1}\psi_1(t,Y_h)h_1(t)dz'|^2_{L^2(D)}&\leq &|\psi_1(t,Y_h)h_1(t)|^2_{L^2(\mathcal{O})}\\
&\leq& C(1+\sup_{t\in[0,T]}|Y_h(t)|^2)|h_1(t)|^2_U\\
&\leq& C|h_1(t)|^2_U,
\end{eqnarray*}
we have
\begin{eqnarray*}
\int_{D}\Delta\bar{v}_h\Big(\int^0_{-1}\psi_1(t,Y_h)h_1(t)dz'\Big)dxdy
&\leq& \varepsilon |\Delta\bar{v}_h|^2_{L^2(D)}+C|h_1(t)|^2_U,
\end{eqnarray*}
thus, we deduce, by (\ref{equ-12}) and (\ref{eq P12 star})
\begin{eqnarray}\label{equ-14}
|\nabla \bar{v}_h(t)|^2+\int^t_0|\Delta \bar{v}_h|^2ds\leq K_4(t),
\end{eqnarray}
where
\[
K_4(t)=e^{K^2_1(t)}\Big[\|v_0\|^2+K_1(t)+K_2(t)+CM\Big].
\]

\textbf{$\Big|\frac{\partial v_h}{\partial z}\Big|^2$ estimates.} \quad Denote $u=\frac{\partial v_h}{\partial z}$. It's clear that $u$ satisfies
\begin{eqnarray}\label{eq-51}
\frac{\partial u}{\partial t}+L_1 u+(v\cdot \nabla)u+\Phi(v)\frac{\partial u}{\partial z}+(u\cdot \nabla)v-(\nabla\cdot v)u+fk\times u-\nabla T=\partial_z (\psi_1(t,Y_h)h_1).
\end{eqnarray}
Taking the inner product of the equation (\ref{eq-51}) with $u$ in $L^2(\mathcal{O})$ and using the boundary condition as Sect. 3.3.2 in Page 13 of \cite{C-T-1}, we get
\begin{eqnarray*}
&&\frac{d|u|^2}{dt}+\frac{3}{2}(|\nabla u|^2+|\partial_z u|^2)\\
&\leq& C(|\nabla \bar{v}_h|^4+|\tilde{v}|^4_6)|u|^2+C|T|^2+|\int_{\mathcal{O}}\partial_z (\psi_1(t,Y_h)h_1)u dxdydz|,
\end{eqnarray*}
 we only need to estimate the following term
\begin{eqnarray*}
\Big|\int_{\mathcal{O}}\partial_z (\psi_1(t,Y_h)h_1)u dxdydz\Big|
&=&\Big|\int_{\mathcal{O}}\psi_1(t,Y_h)h_1\partial_z u dxdydz\Big|\\
&\leq& C|\psi_1(t,Y_h)h_1(t)|_{L^2(\mathcal{O})}|\partial_zu|\\
&\leq& \varepsilon|\partial_zu|^2+C|\psi_1(t,Y_h)h_1(t)|^2_{L^2(\mathcal{O})}.
\end{eqnarray*}
By \textbf{Hypothesis H0} and (\ref{equ-12}) again
\begin{eqnarray}\label{eq P13 star}
|\psi_1(t,Y_h)h_1(t)|^2\leq C|h_1(t)|^2_U,
\end{eqnarray}
thus, similar as equation (75) in \cite{C-T-1},
\begin{eqnarray}\label{equ-15}
\Big|\frac{\partial v_h}{\partial z}\Big|^2+\int^t_0\Big|\nabla\frac{\partial v_h}{\partial z} \Big|^2ds+\int^t_0\Big|\frac{\partial^2 v_h}{\partial z^2}\Big |^2ds\leq K_5(t),
\end{eqnarray}
where
\[
K_5(t)=e^{\Big(K^2_4(t)+K^{\frac{3}{2}}_3(t)\Big)t}\Big[ \|v_0\|^2+K_1(t)+CM\Big].
\]
\textbf{$|\nabla v_h|^2$ estimates.} \quad Taking the inner product of the equation (\ref{eq-47}) with $-\Delta v_h$ in $L^2(\mathcal{O})$. As Page 14 in \cite{C-T-1}, we reach
\begin{eqnarray*}
&&\frac{d|\nabla v_h|^2}{dt}+\frac{3}{2}(|\Delta v_h|^2+|\nabla \partial_z v_h|^2)\\
&\leq& C(|v|^4_6+|\nabla v_h|^2|\partial_z v_h|^2)|\nabla v_h|^2+C|\nabla T_h|^2+|\int_{\mathcal{O}}\psi_1(t,Y_h)h_1(t)\Delta v_h dxdydz|,
\end{eqnarray*}
we only need to estimate the following additional term
\begin{eqnarray*}
|\int_{\mathcal{O}}\psi_1(t,Y_h)h_1(t)\Delta v_h dxdydz|
&\leq& C|\psi_1(t,Y_h)h_1(t)|_{L^2(\mathcal{O})}|\Delta v_h|\\
&\leq& \varepsilon |\Delta v_h|^2+C|\psi_1(t,Y_h)h_1(t)|^2.
\end{eqnarray*}
By (\ref{eq P13 star}), as equation (77) in \cite{C-T-1},
\begin{eqnarray}\label{equ-16}
|\nabla v_h(t)|^2+\int^t_0|\Delta v_h(s)|^2ds+\int^t_0|\nabla \frac{\partial v_h(t)}{\partial z}|^2dt\leq K_6(t),
\end{eqnarray}
where
\[
K_6(t)=e^{\Big(K^{\frac{2}{3}}_4(t)t+K_1(t)K_5(t)\Big)}\Big[\|v_0\|^2+K_1(t)+CM\Big].
\]
\textbf{$\|T\|$ estimates.} \quad Taking the inner product of the equation (\ref{eq-48}) with $-\Delta T-T_{zz}$ in $L^2(\mathcal{O})$, in the same way as Sect. 3.3.4  in  Page 15 of \cite{C-T-1}, we get
\begin{eqnarray*}
&&\frac{d(|\nabla T_h|^2+|\partial_z T_h|^2)}{dt}+\frac{3}{2}(|\Delta T_h|^2+|\nabla \partial_z T_h|^2+|\partial_{zz} T_h|^2)\\
&\leq& C(|v_h|^4_6+|\nabla v_h|^2|\Delta v_h|^2)(|\nabla T_h|^2+|\partial_z T_h|^2)+|\int_{\mathcal{O}}\psi_2(t,Y_h)h_2(t)(\Delta T_h+\partial^2_z T_h )dxdydz|,
\end{eqnarray*}
we only need to estimate the following term
\begin{eqnarray*}
&&|\int_{\mathcal{O}}\psi_2(t,Y_h)h_2(t)(\Delta T_h+\partial^2_z T_h )dxdydz|\\
&\leq& C|\psi_2(t,Y_h)h_2(t)|_{L^2(\mathcal{O})}|\Delta T_h+\partial^2_z T_h |\\
&\leq & \varepsilon(|\Delta T_h|^2+|\partial^2_z T_h|^2+|\nabla \partial_z T_h|^2)+C|h_2(t)|^2_U,
\end{eqnarray*}
where \textbf{Hypothesis H0} is used.
Thus, we obtain
\begin{eqnarray}\label{eq-52}
|\nabla T_h|^2+|\partial_z T_h|^2+\int^t_0\Big(|\Delta T_h|^2+|\nabla \partial_z T_h|^2+|\partial_{zz} T_h|^2\Big)dt\leq K_7(t),
\end{eqnarray}
where
\[
K_7(t)=e^{\Big(K^{2}_4(t)t+K^2_6(t)\Big)}\Big[\|T_0\|^2+CM\Big].
\]

\subsubsection{Proof of Theorem \ref{thm-4}}
Now, we are ready to prove Theorem \ref{thm-4}.

\textbf{Proof of Theorem \ref{thm-4}} \quad
Combining (\ref{equ-13})-(\ref{eq-52}) and using proof by contradiction, we obtain the global existence of strong solution of (\ref{equ-9}).

 In the following, we only need to prove the uniqueness and continuously dependence on the initial data. Let $Y_1=(v_1,T_1,p^1_b), Y_2=(v_2,T_2,p^2_b)$ be two strong solutions of (\ref{equ-9}), for convenience, here, we omit the index $h$. Denote $r=v_1-v_2, \eta=T_1-T_2, q_b=p^1_b-p^2_b$, it is clear that
\begin{eqnarray}\label{eq-53}
\frac{dr}{d t}&+&L_1r+(v_1\cdot\nabla)r+(r\cdot \nabla) v_2+\Phi(v_1)\frac{\partial r}{\partial z}+\Phi(r)\frac{\partial v_2}{\partial z}+fk\times r\\ \notag
&&+\nabla q_b-\int^z_{-1}\nabla\eta(x,y,z',t)dz'=\psi_1(t,Y_1(t))h_1-\psi_1(t,Y_2(t))h_1,\\
\label{eq-54}
\frac{d\eta}{d t}&+&L_1\eta+(v_1\cdot\nabla)\eta+(r\cdot \nabla) T_2+\Phi(v_1)\frac{\partial \eta}{\partial z}+\Phi(r)\frac{\partial T_2}{\partial z}\\ \notag
&&=\psi_2(t,Y_1(t))h_2(t)-\psi_2(t,Y_2(t))h_2,\\
\label{eq-55}
&&r(x,y,z,0)=v^1_0-v^2_0,\\
\label{eq-56}
&&\eta(x,y,z,0)=T^1_0-T^2_0.
\end{eqnarray}
\textbf{$L^2$ estimates of $r$.} \quad Taking the inner product of the equation (\ref{eq-53}) with $r$ in $L^2(\mathcal{O})$ as Page 16 in \cite{C-T-1}, we get
\begin{eqnarray*}
&&\frac{d|r|^2}{dt}+\frac{3}{2}(|\nabla r|^2+|\partial_z r|^2)\\
&\leq& C|\nabla v_2|^4|r|^2+C|r|^2|\partial_z v|^2|\nabla v_z|^2+|\int_{\mathcal{O}}\Big(\psi_1(t,Y_1(t))h_1(t)-\psi_1(t,Y_2(t))h_1(t)\Big)rdxdydz|,
\end{eqnarray*}
we only need to estimate the additional term,
\begin{eqnarray*}
&&|\int_{\mathcal{O}}\Big(\psi_1(t,Y_1(t))h_1(t)-\psi_1(t,Y_2(t))h_1(t)\Big)rdxdydz|\\
&\leq & |\psi_1(t,Y_1(t))h_1(t)-\psi_1(t,Y_2(t))h_1(t)||r|\\
&\leq& |\psi_1(t,Y_1(t))-\psi_1(t,Y_2(t))|_{\mathcal{L}_2(U,H)}|h_1(t)|_U|r|\\
&\leq& |Y_1(t)-Y_2(t)||h_1(t)|_U|r|\\
&\leq& |h_1(t)|_U|r|^2+|h_1(t)|_U|r||\eta|\\
&\leq& C(1+|h_1(t)|^2_U)|r|^2+C|h_1(t)|^2_U|\eta|^2,
\end{eqnarray*}
Similarly, we can obtain $L^2$ estimate of $\eta$. \quad
Taking the inner product of the equation (\ref{eq-54}) with $\eta$ in $L^2(\mathcal{O})$, we reach
\begin{eqnarray*}
&&\frac{d|\eta|^2}{dt}+\frac{3}{2}(|\nabla \eta|^2+|\partial_z \eta|^2)\\
&\leq& C|\nabla T_2|^4|\eta|^2+C|\eta|^2|\partial_z T_2|^2|\nabla \partial_z T_2|^2+|\int_{\mathcal{O}}\Big(\psi_1(t,Y_1(t))h_1(t)-\psi_1(t,Y_2(t))h_2(t)\Big)\eta dxdydz|,
\end{eqnarray*}
we only need to estimate the additional term,
\begin{eqnarray*}
&&|\int_{\mathcal{O}}\Big(\psi_1(t,Y_1(t))h_2(t)-\psi_1(t,Y_2(t))h_2(t)\Big)\eta dxdydz|\\
&\leq & |\psi_2(t,Y_1(t))h_2(t)-\psi_2(t,Y_2(t))h_2(t)||\eta|\\
&\leq& \|\psi_2(t,Y_1(t))-\psi_2(t,Y_2(t))\|_{\mathcal{L}_2(U,H)}|h_2(t)|_U|\eta|\\
&\leq& |Y_1(t)-Y_2(t)||h_2(t)|_U|\eta|\\
&\leq& C(1+|h_2(t)|^2_U)|\eta|^2+C|h_2(t)|^2_U|r|^2,
\end{eqnarray*}
therefore, we have
\[
|r(t)|^2+|\eta(t)|^2\leq(|r(0)|^2+|\eta(0)|^2)e^{C\Big(K^2_6 t+K^2_7 t+K_5K_6+K^2_6+C(1+M)\Big)}.
\]
The above inequality proves the continuous dependence of the solutions on the initial data, and in particular, when $r(0)=\eta(0)=0$, we have $r(t)=\eta(t)=0$, for all $t\geq0$. Therefore, the strong solution is unique.

$\hfill\blacksquare$

(\ref{equ-16}) and (\ref{eq-52}) imply
\begin{cor}\label{cor-2}
Let $Y_{h}$ be the unique strong solution of (\ref{equ-9}) with $h\in T_M$, then
\[
\sup_{0\leq t\leq T}\|Y_{h}(t)\|^2+\int^T_0\|Y_{h}(t)\|^2_2dt\leq C(T,M,\|Y_0\|).
\]
\end{cor}
Now, define $\mathcal{G}^0: C([0,T];U)\rightarrow \Re$ by
\begin{eqnarray}\label{equ-28}
\mathcal{G}^0(\tilde{h})=\left\{
                 \begin{array}{ll}
                   Y_h, & \rm{if}\  \tilde{h}=\int^{\cdot}_{0}h(s)ds \ \rm{for}\ \rm{some}\ h\in L^2([0,T]; U), \\
                   0, & \rm{otherwise}.
                 \end{array}
               \right.
\end{eqnarray}

\subsection{Compactness of $Y_n$ }
Let $Y_{n}$ be the unique strong solution of (\ref{equ-9}) with $h_n\in T_M$ and $h_n=(h^1_n,h^2_n)$, in this section, we aim to prove $(Y_n)_{n\in \mathbb{N}_+}$ is compact in $L^2([0,T];V)$.  Refer to \cite{FG95}, we introduce the following definition and lemma which are needed below.

Given $p>1, \alpha\in (0,1)$, let $W^{\alpha,p}([0,T]; K)$ be the Sobolev space of all $u\in L^p(0,T;K)$ such that
\[
\int^T_0\int^T_0\frac{|u(t)-u(s)|_K^{ p}}{|t-s|^{1+\alpha p}}dtds< \infty,
\]
endowed with the norm
\[
|u|^p_{W^{\alpha,p}(0,T; H)}=\int^T_0|u(t)|_K^pdt+\int^T_0\int^T_0\frac{|u(t)-u(s)|_K^{ p}}{|t-s|^{1+\alpha p}}dtds.
\]
\begin{lemma}\label{lem-1}
Let $B_0\subset B\subset B_1$ be Banach spaces, $B_0$ and $B_1$ reflexive, with compact embedding of $B_0$ in $B$. Let $p\in (1, \infty)$ and $\alpha \in (0, 1)$ be given. Let $X$ be the space
\[
X= L^p([0, T]; B_0)\cap W^{\alpha, p}([0,T]; B_1),
\]
endowed with the natural norm. Then the embedding of $X$ in $L^p([0,T];B)$ is compact.
\end{lemma}
Now, we will apply Lemma \ref{lem-1} to obtain compactness of $(Y_n)_{n\in \mathbb{N}_+}$.
\begin{prp}\label{prp-1}
$(Y_n)_{n\geq 0}$ are compact in $L^2([0,T];V)$.
\end{prp}
\textbf{Proof of Proposition \ref{prp-1}} \quad From (\ref{equ-9}), we have
\begin{eqnarray*}
Y_n(t)&=&Y_0-\int^t_0AY_n(s)ds-\int^t_0B(Y_n(s),Y_n(s))ds
-\int^t_0G(Y_n(s))ds+\int^t_0\psi(s,Y_n(s))h_n(s)ds
\\  &=& J^1_n+J^2_n(t)+J^3_n(t)+J^4_n(t)+J^5_n(t).
\end{eqnarray*}
Refer to Sect. 4.2 in \cite{RR}, we have
\begin{eqnarray*}
|J^1_n|^2&\leq& C_1, \\
|J^4_n|^2_{W^{\alpha,2}(0,T; V')}&\leq& C\left(\sup_{0\leq s\leq T}|Y_n(s)|^2\right)\leq C_{2,\alpha},\quad \alpha\in (0,\frac{1}{2}),\\
|J^2_n|^2_{W^{\alpha,2}(0,T; V')}&\leq& C \int^{T}_{0}\|Y_n(s)\|^2ds\leq C_{3,\alpha}\quad \alpha\in (0,\frac{1}{2}),
\end{eqnarray*}
for suitable positive constants $C_1$, $C_2$, $C_3$.
For $J^3_n$, by Lemma \ref{le-2},
\[
\|B(Y,Y_1)\|_{V'}\leq C\|Y\|\|Y_1\|_2,
\]
then
\begin{eqnarray}\label{eq Page 17 1}
|B(Y_n,Y_n)|^2_{L^2(0,T;V')}\leq C_4\Big(\sup_{0\leq s\leq T}\|Y_n(s)\|^2\Big)\int^{T}_{0}\|Y_n(s)\|^2_2ds,
\end{eqnarray}
by Corollary \ref{cor-2}, we obtain
\[
|J^3_n|^2_{W^{\alpha,2}(0,T; V')}\leq C_{5,\alpha} \quad \alpha\in (0,1).
\]
As to $J^5_n$, since
\begin{eqnarray}\label{equ-25}
|\int^t_s\psi(u,Y_n(u)) h_n(u)du|^2
&\leq& \int^t_s|\psi(u,Y_n(u)) h_n(u)|^2du\\ \notag
&\leq& \int^t_s|h_n(u)|^2_U(1+|Y_n(u)|^2)du\\ \notag
&\leq &\int^t_s|h_n(u)|^2_Udu+\int^t_s|h_n(u)|^2_U|Y_n(u)|^2du\\ \notag
&\leq &\Big(1+\sup_{0\leq u\leq T}|Y_n(u)|^2\Big)\int^t_s|h_n(u)|^2_Udu\\ \notag
&\leq& C\int^t_s|h_n(u)|^2_Udu,
\end{eqnarray}
then, by Fubini Theorem,
\begin{eqnarray*}
|J^5_n|^2_{W^{\alpha,2}(0,T; H)}&=&|\int^t_0\psi(s,Y_n(s)) h_n(s)ds|^2_{W^{\alpha,2}(0,T; H)}\\
&=&\int^T_0|\int^t_0\psi(Y_n) h_n(s)ds|^2dt +\int^T_0\int^T_0\frac{|\int^t_s\psi(u,Y_n(u)) h_n(u)du|^2}{|t-s|^{1+2\alpha}}dtds\\
&\leq& C_6(\alpha,M),\ \ \forall \alpha\in (0,\frac{1}{2}).
\end{eqnarray*}

Collecting all the previous inequalities we obtain
\[
|Y_n|^2_{W^{\alpha,2}([0,T];V')}\leq C_7(\alpha),\ \ \  \forall \alpha\in (0,\frac{1}{2})
\]
for some constant $C_7(\alpha)>0$. Recalling Corollary \ref{cor-2}, we have that $Y_n$ are bounded uniformly in $n$ in the space
 \[
 L^2([0,T];D(A))\cap W^{\alpha,2}([0,T];V'),
 \]
 by Lemma \ref{lem-1}, $Y_n$ are compact in $L^2([0,T];V)$. $\hfill\blacksquare$

\begin{cor}\label{cor-3}
There exists a subsequence still denoted by $Y_n$ and $\check{Y}\in L^{\infty}([0,T];V)\cap L^2([0,T];V)\cap L^2([0,T];D(A))$ such that
\begin{eqnarray*}
&Y_n&\rightarrow \check{Y} \ weakly\  star \ in\ L^{\infty}([0,T];V),\\
&Y_n&\rightarrow \check{Y} \ strongly \ in\ L^2([0,T];V),\\
&Y_n&\rightharpoonup\check{Y} \ weakly \ in\ L^2([0,T];D(A)).
\end{eqnarray*}
\end{cor}
\subsection{The Property of $\check{Y}$}
Fix a sequence $(h_n)_{n\geq 0}$ such that $h_n\rightharpoonup h$ weakly in $T_M$, from Theorem \ref{thm-4} and Corollary \ref{cor-3}, the limit of $Y_n$ exists and we denote it by $\check{Y}$. The following proposition tells us that $\check{Y}$ is the solution of (\ref{equ-9}) with $h$.
\begin{prp}\label{prp-2}
The above $\check{Y}$ satisfies
\begin{eqnarray}\label{equ-23}
\left\{
  \begin{array}{ll}
    d\check{Y}(t)+A\check{Y}(t)dt+B(\check{Y}(t),\check{Y}(t))dt+G(\check{Y}(t))dt=\psi(t,\check{Y}(t))h(t)dt,\\
   \check{Y}(0)=Y_0,
  \end{array}
\right.
\end{eqnarray}
\end{prp}
Before the proof, we firstly give a lemma for the nonlinear term.
\begin{lemma}\label{lem-2}
Let $w\in D(A^{\frac{3}{2}})$, $u_{\nu}\rightarrow u$ strongly in $L^2([0,T];V)$ as $\nu\rightarrow 0$, then
\[
\int^T_0(B(u_{\nu}(t),u_{\nu}(t)),w(t))dt\rightarrow \int^T_0(B(u(t),u(t)),w(t))dt, \ as\ \nu\rightarrow 0.
\]
\end{lemma}
\textbf{Proof of Lemma \ref{lem-2}} \quad
\begin{eqnarray*}
&&\Big|\int^T_0(B(u_{\nu}(t),u_{\nu}(t)),w)dt-\int^T_0(B(u(t),u(t)),w)dt\Big|\\
&&\leq \int^T_0|(B(u_{\nu},u_{\nu}-u),w)|dt+\int^T_0|(B(u_{\nu}-u,u),w)|dt\\
&&:=I_1+I_2,
\end{eqnarray*}
refer to \cite{STW},
\[
\|B(Y,Y_1)\|_{-3}\leq C|Y|\|Y_1\|,
\]
then, by H\"{o}lder inequality and Sobolev embedding, we have
\begin{eqnarray*}
I_1&\leq& C\int^T_0\|u_{\nu}\||u_{\nu}-u||w|_{D(A^{\frac{3}{2}})}dt\\
&\leq &C|w|_{D(A^{\frac{3}{2}})}\int^T_0\|u_{\nu}\||u_{\nu}-u|dt\\
&\leq&C|w|_{D(A^{\frac{3}{2}})}\Big(\int^T_0\|u_{\nu}\|^2dt\Big)^{\frac{1}{2}}\Big(\int^T_0|u_{\nu}-u|^2dt\Big)^{\frac{1}{2}},\\
I_2&\leq& C\int^T_0\|u_{\nu}-u\||u||w|_{D(A^{\frac{3}{2}})}dt\\
&\leq &C|w|_{D(A^{\frac{3}{2}})}\int^T_0\|u_{\nu}-u\||u|dt\\
&\leq & C|w|_{D(A^{\frac{3}{2}})}\Big(\int^T_0\|u_{\nu}-u\|^2dt\Big)^{\frac{1}{2}}\Big(\int^T_0|u|^2dt\Big)^{\frac{1}{2}},
\end{eqnarray*}
since $u_{\nu}\rightarrow u$ strongly in $L^2([0,T];V)$, we have $I_1+I_2\rightarrow 0$.

$\hfill\blacksquare$

\textbf{Proof of Proposition \ref{prp-2}} \quad
Denoting an orthonormal basis of $D(A^{\frac{3}{2}}_1)$ by $\{w^1_j\}_{j\geq 1}$ and an orthonormal basis of $D(A^{\frac{3}{2}}_2)$ by $\{w^2_j\}_{j\geq 1}$, then following Sect. \ref{sec 2}, we obtain an orthonormal basis of $D(A^{\frac{3}{2}})$, which is denoted by $\{w_j\}_{j\geq 1}$. Taking a test function $\phi(t)$ a continuously differentiable on $[0,T]$ satisfying $\phi(T)=0$. From (\ref{equ-9}), we have
\begin{eqnarray*}
&\int^T_0(\frac{dY_n}{dt},\phi(t)w_j)dt+\int^T_0(AY_n,\phi(t)w_j)dt+\int^T_0(B(Y_n,Y_n),\phi(t)w_j)dt\\
&+\int^T_0(G(Y_n),\phi(t)w_j)dt
=\int^T_0(\psi(t,Y_n(t))h_n(t),\phi(t)w_j)dt,
\end{eqnarray*}
by integration by parts,
\begin{eqnarray*}
&-(Y_0,\phi(0)w_j)-\int^T_0(Y_n(t),\phi'(t)w_j)dt+\int^T_0(Y_n(t),\phi(t)Aw_j)dt+\int^T_0(B(Y_n,Y_n),\phi(t)w_j)dt\\
&+\int^T_0(G(Y_n),\phi(t)w_j)dt=\int^T_0(\psi(t,Y_n(t)) h_n(t),\phi(t)w_j)dt.
\end{eqnarray*}
Denote the above equality by symbols that $J_1+J_2+J_3+J_4+J_5=J_6$, in the following, we will estimate these terms one by one.

For $J_2+J_3$, by H\"{o}lder inequality and  $Y_n \rightarrow \check{Y}$ strongly in $L^2([0,T];V)$, we have
\[
J_2+J_3 \rightarrow -\int^T_0(\check{Y}(t),\phi'(t)w_j)dt+\int^T_0(\check{Y}(t),\phi(t)Aw_j)dt.
\]
For $J_4$, it follows from Lemma \ref{lem-2} that
 \[
 J_4\rightarrow\int^T_0(B(\check{Y},\check{Y}),\phi(t)w_j)dt, \ n\rightarrow \infty.
 \]
 As to $J_5$, denote $\check{Y}=(\check{v},\check{T})$, we have
\begin{eqnarray*}
&&\int^T_0(G(Y_n),\phi(t)w_j)dt-\int^T_0(G(\check{Y}),\phi(t)w_j)dt\\
&=&\int^T_0\Big(fk\times(v_n-\check{v}),\phi(t)w^1_j\Big)dt
+\int^T_0\Big(\int^z_{-1}\nabla(T_n-\check{T})dz',\phi(t)w^1_j\Big)dt\\
&=&K_1+K_2.
\end{eqnarray*}
For $K_1$, by H\"{o}lder inequality and Corollary \ref{cor-3}, we have $K_1\rightarrow 0$. For $K_2$, we have
 \begin{eqnarray*}
\int^T_0\Big(\int^z_{-1}\nabla(T_n-\check{T})dz',\phi(t)w^1_j\Big)dt
&=&-\int^T_0\Big(\int^z_{-1}(T_n-\check{T})dz',\phi(t)\nabla w^1_j\Big)dt\\
&\leq &\int^T_0|T_n-\check{T}||\phi(t)\nabla w^1_j|dt\rightarrow 0,
\end{eqnarray*}
thus,
\[
J_5\rightarrow \int^T_0(G(\check{Y}),\phi(t)w_j)dt, \ n\rightarrow \infty.
\]
For $J_6$,
\begin{eqnarray*}
&&\Big|\int^T_0\Big(\psi(t,Y_n)h_n(t),\phi(t)w_j\Big)dt-\int^T_0\Big(\psi(t,\check{Y}))h(t),\phi(t)w_j\Big)dt\Big|\\ &&\quad \leq\Big|\int^T_0\Big((\psi(t,Y_n)-\psi(t,\check{Y}))h_n(t),\phi(t)w_j\Big)dt\Big|
+\Big|\int^T_0\Big(\psi(t,\check{Y})(h_n(t)-h(t)),\phi(t)w_j\Big)dt\Big|,\\
&&:= K_3+K_4.
\end{eqnarray*}
By H\"{o}lder inequality, for $K_3$,we have
\begin{eqnarray*}
 \Big|\int^T_0\Big((\psi(t,Y_n)-\psi(t,\check{Y}))h_n(t),\phi(t)w_j\Big)dt\Big|
 &\leq& \int^T_0|(\psi(t,Y_n)-\psi(t,\check{Y}))h_n(t)||\phi(t)w_j|dt\\
 &\leq& \int^T_0|\psi(t,Y_n)-\psi(t,\check{Y}))|_{\mathcal{L}_2(U,H)}|h_n(t)|_U|\phi(t)w_j|dt\\
 &\leq& C \int^T_0|Y_n-\check{Y}||h_n(t)|_Udt\\
 &\leq& C(\int^T_0|Y_n-\check{Y}|^2dt)^{\frac{1}{2}}(\int^T_0|h_n(t)|_U^2dt)^{\frac{1}{2}},
\end{eqnarray*}
since $Y_n\rightarrow \check{Y}$ strongly in $L^2([0,T];V)$, we have
$K_3\rightarrow 0, n\rightarrow \infty$.


$h_n-h\rightarrow 0$ weakly in $L^2([0,T];U)$ and Corollary \ref{cor-3} imply that $K_4\rightarrow 0$, $n\rightarrow \infty$.
Thus,
\[
J_6\rightarrow\int^T_0\Big(\psi(t,\check{Y}))h(t),\phi(t)w_j\Big)dt, \ n\rightarrow \infty.
\]
Hence, we have
\begin{eqnarray}\label{equ-24}
&-\int^T_0(\check{Y}(t),\phi'(t)w_j)dt+\int^T_0(\check{Y}(t),Aw_j\phi(t))dt
+\int^T_0(B(\check{Y},\check{Y}),\phi(t)w_j)dt
+\int^T_0(G(\check{Y}),\phi(t)w_j)dt\\
&=(Y_0,\phi(0)w_j)+\int^T_0(\psi(t,\check{Y}(t)) h(t),\phi(t)w_j)dt. \notag
\end{eqnarray}
Since the above equality holds for each $j$, so (\ref{equ-24}) holds for any $\zeta$, which is a finite linear combination of $w_j$, that is
\begin{eqnarray}\label{eq-80}
&-\int^T_0(\check{Y}(t),\phi'(t)\zeta)dt+\int^T_0(\check{Y}(t),A\phi(t)\zeta)dt
+\int^T_0(B(\check{Y},\check{Y}),\phi(t)\zeta)dt
+\int^T_0(G(\check{Y}),\phi(t)\zeta)dt\\
&\quad\quad\quad=({Y}_0,\phi(0)\zeta)+\int^T_0(\psi(t,\check{Y}(t)) h(t),\phi(t)\zeta)dt. \notag
\end{eqnarray}
Since $D(A^{\frac{3}{2}})$ is dense in $V$, we have the following equality holds as an equality in the distribution sense in $L^2([0,T];V')$,
\begin{eqnarray}\label{eq-81}
\frac{d}{dt}(\check{Y},\zeta)+(A\check{Y},\zeta)+(B(\check{Y},\check{Y}),\zeta)+(G(\check{Y}),\zeta)=(\psi(t,\check{Y}),\zeta),
\end{eqnarray}
which is exactly (\ref{equ-23}).
 Finally, it remains to prove $\check{Y}(0)={Y}_0$. For this, multiplying (\ref{eq-81}) with the same $\phi(t)$ as above, integrate with respect to $t$, and integrate by parts, we have
\begin{eqnarray}\label{eq-79}
&-\int^T_0(\check{Y}(t),\phi'(t)\zeta)dt+\int^T_0(\check{Y}(t),A\phi(t)\zeta)dt
+\int^T_0(B(\check{Y},\check{Y}),\phi(t)\zeta)dt
+\int^T_0(G(\check{Y}),\phi(t)\zeta)dt
\\
&\quad=(\check{Y}(0),\phi(0)\zeta)+\int^T_0(\psi(t,\check{Y}(t)) h(t),\phi(t)\zeta)dt\notag.
\end{eqnarray}
 By comparison with (\ref{eq-80}), we see that  $(\check{Y}(0)-{Y}_0,\phi(0)\zeta)=0$ for each $\zeta\in D(A^{\frac{3}{2}})$ and for each function $\phi$ of the type considered. We can choose $\phi$ such that $\phi(0)\neq 0$, therefore,
\[
(\check{Y}(0)-{Y}_0,\zeta)=0,\  \forall \zeta\in D(A^{\frac{3}{2}}).
\]
As $D(A^{\frac{3}{2}})$ is dense in $V$, we have that $\check{Y}(0)={Y}_0$, which conclude the result.
$\hfill\blacksquare$

\subsection{The Continuity of $\check{Y}$ in V}
In this section, we will use the following Lemma \ref{lem-4} to obtain  Proposition \ref{prp-4}.
\begin{lemma}\label{lem-4}
For $V$ and $H$ are Hilbert spaces ($V'$ is the dual space of $V$) with $V\subset\subset H=H'\subset V'$, where $V\subset\subset H$ denotes $V$ is compactly embedded in $H$. If $u\in L^2([0,T];V)$, $\frac{du}{dt}\in L^2([0,T];V')$, then $u\in C([0,T];H)$.
\end{lemma}
\begin{prp}\label{prp-4}
$\check{Y}\in C([0,T];V)$.
\end{prp}

\textbf{Proof of Proposition \ref{prp-4}} \quad
 Following Lemma \ref{lem-4}, we should firstly prove that $\frac{d\check{Y}}{dt}$ is in $L^2([0,T];V')$. Indeed,
 in the proof of Proposition \ref{prp-2}, we know $\check{Y}\in L^2([0,T],D(A))\cap L^{\infty}([0,T],V)$ and
\begin{eqnarray*}
\frac{d\check{Y}}{dt}=-A\check{Y}-B(\check{Y},\check{Y})-G(\check{Y})+\psi(t,\check{Y})h.
\end{eqnarray*}
For $A\check{Y}$, since $\check{Y}$ is bounded in $L^2([0,T];D(A))$ and $A$ is continuous linear operator from $D(A)$ to $H$, thus, $A\check{Y}$ is bounded in $L^2([0,T];H)$.
For $B(\check{Y},\check{Y})$,  similar to (\ref{eq Page 17 1}), we have
\[
\|B(\check{Y},\check{Y})\|_{L^2([0,T];V')}\leq C.
\]
For $G(\check{Y})$, we have
\[
|G(\check{Y})|^2_{L^2([0,T]; V')}\leq C\left(\sup_{0\leq s\leq T}|\check{Y}(s)|^2\right)\leq C.
\]
 For $\psi(t,\check{Y})h$, by \textbf{Hypothesis H0}, we have
 \begin{eqnarray*}
 |\psi(t,\check{Y})h|^2_{L^2([0,T];H)}
 &=& \int^T_0 |\psi(t,\check{Y})h|^2 dt\\
 &\leq& \int^T_0 \|\psi(t,\check{Y})\|^2_{\mathcal{L}_2(U;H)}|h|^2_U dt\\
 &\leq&  \int^T_0 (1+|\check{Y}|^2)|h|^2_U dt\\
 &\leq&  C\sup_{t\in [0,T]} |\check{Y}(t)|^2\int^T_0 |h|^2_U dt
 \leq  CM.
 \end{eqnarray*}
  Collecting all the above estimates, we get
  \[
   \frac{d\check{Y}}{dt}\ \in \ L^2([0,T];V').
  \]
  Recalling Corollary \ref{cor-2}, $\check{Y} \in L^2([0,T];D(A))$ and applying Lemma \ref{lem-4}, we conclude the result.
  $\hfill\blacksquare$

By the uniqueness of (\ref{equ-9}), we have the following corollary.
\begin{cor}\label{cor-1}
 $\check{Y}=Y_h$, where $Y_h$ is the unique strong solution of (\ref{equ-9}) with $h$.
\end{cor}
Moreover, we have
\begin{thm}\label{thm-6}
$Y_n-\check{Y}\rightarrow 0$ in $\Re$ as $n\rightarrow \infty$.
\end{thm}
\textbf{Proof of Theorem \ref{thm-6}} \quad Denote $Y_n=(v_n, T_n, p_n)$ with $h_n=(h^1_n,h^2_n)$ and $\check{Y}=(\check{v},\check{T},\check{p}_b)$ with $h=(h_1,h_2)$. Let
$r_n=v_n-\check{v}, \eta_n=T_n-\check{T}, q_n=p_n-\check{p}_b$, then
\begin{eqnarray}\label{equ-29}
\frac{dr_n}{d t}&+&A_1r_n+(v_n\cdot\nabla)r_n+(r_n\cdot \nabla) \check{v}+\Phi(v_n)\frac{\partial r_n}{\partial z}+\Phi(r_n)\frac{\partial \check{v}}{\partial z}+fk\times r_n\\ \notag
&&+\nabla q_n-\int^z_{-1}\nabla\eta_n(x,y,z',t)dz'=\psi_1(t,Y_n(t))h^1_n-\psi_1(t,\check{Y})h_1,\\
\label{equ-30}
\frac{d\eta_n}{d t}&+&A_2\eta_n+(v_n\cdot\nabla)\eta_n+(r_n\cdot \nabla) \check{T}+\Phi(v_n)\frac{\partial \eta_n}{\partial z}+\Phi(r_n)\frac{\partial \check{T}}{\partial z}\\ \notag
&&=\psi_2(t,Y_n)h^2_n-\psi_2(t,\check{Y})h_2,\\
\label{equ-31}
&&r_n(x,y,z,0)=0,\\
\label{equ-32}
&&\eta_n(x,y,z,0)=0.
\end{eqnarray}
\textbf{$H^1$ estimate of $r_n$}. \quad Taking the inner product of (\ref{equ-29}) with $A_1 r_n$ in $L^2(\mathcal{O})$, then integrating the time from 0 to $t$, it reaches

\begin{eqnarray*}
&&\|r_n(t)\|^2+2\int^t_0\|r_n(s)\|^2_2ds\\
&=&-2\int^t_0\Big((v_n\cdot\nabla)r_n+\Phi(v_n)\frac{\partial r_n}{\partial z}, A_1r_n\Big)ds\\
&&-2\int^t_0\Big((r_n\cdot \nabla) \check{v}+\Phi(r_n)\frac{\partial \check{v}}{\partial z},A_1 r_n\Big)ds\\
&&-2\int^t_0\Big(fk\times r_n+\nabla q_n-\int^z_{-1}\nabla\eta_n(x,y,z',t)dz',A_1r_n\Big)ds\\
&&+2\int^t_0\Big(\psi_1(s,Y_n)h^1_n-\psi_1(s,\check{Y})h_1,A_1 r_n\Big)ds\\
&:=&I_1+I_2+I_3+I_4.
\end{eqnarray*}
Applying H\"{o}lder inequality, Lemma \ref{le-2} and Corollary \ref{cor-2} to $I_1$ and $I_2$, we obtain
\begin{eqnarray*}
&& \Big|\int^t_0\Big((v_n\cdot\nabla)r_n+\Phi(v_n)\frac{\partial r_n}{\partial z}, A_1r_n\Big)ds\Big|\\
&\leq& C\int^t_0\|r_n(s)\|_2|\nabla r_n||v_n|_{\infty}ds +C\int^t_0\|r_n(s)\|_2\|v_n\|^{\frac{1}{2}}\|v_n\|^{\frac{1}{2}}_2|\frac{\partial r_n}{\partial z}|^{\frac{1}{2}}|\nabla \frac{\partial r_n}{\partial z}|^{\frac{1}{2}}ds\\
&\leq& \varepsilon\int^t_0\|r_n(s)\|^2_2ds+C\int^t_0(1+\|v_n\|^2)\|v_n\|^2_2\|r_n\|^2ds,
\end{eqnarray*}
and
\begin{eqnarray*}
&& \Big|\int^t_0\Big((r_n\cdot \nabla) \check{v}+\Phi(r_n)\frac{\partial \check{v}}{\partial z},A_1 r_n\Big)ds\Big|\\
&\leq& C\int^t_0\|r_n(s)\|_2\|\check{v} \|^{\frac{1}{2}}\|\check{v} \|^{\frac{1}{2}}_2\|r_n\|ds +C\int^t_0\|r_n(s)\|_2\|r_n\|^{\frac{1}{2}}\|r_n\|^{\frac{1}{2}}_2|\frac{\partial \check{v} }{\partial z}|^{\frac{1}{2}}|\nabla \frac{\partial \check{v}}{\partial z}|^{\frac{1}{2}}ds\\
&\leq& \varepsilon\int^t_0\|r_n(s)\|^2_2ds+C\int^t_0(1+\|\check{v}\|^2)\|\check{v}\|^2_2\|r_n\|^2ds.
\end{eqnarray*}
For $I_3$,
\begin{eqnarray*}
&& \Big|\int^t_0\Big(fk\times r_n+\nabla q_n-\int^z_{-1}\nabla\eta_n(x,y,z',t)dz',A_1r_n\Big)ds\Big|\\
&\leq& \varepsilon\int^t_0\|r_n(s)\|^2_2ds+C\int^t_0(|r_n|^2+|\nabla \eta_n|^2)ds.
\end{eqnarray*}
Finally, for $I_4$,
\begin{eqnarray*}
&& \int^t_0\Big(\psi_1(s,Y_n)h^1_n-\psi_1(s,\check{Y})h_1,A_1 r_n\Big)ds\\
&=& \int^t_0\Big((\psi_1(s,Y_n)-\psi_1(s,\check{Y}))h^1_n, A_1 r_n\Big)ds +\int^t_0\Big(\psi_1(s,\check{Y})(h^1_n-h_1), A_1 r_n\Big)ds\\
&:=& J_1+J_2,
\end{eqnarray*}
by H\"{o}lder inequality, the Young inequality and \textbf{Hypothesis H0}, we have
\begin{eqnarray*}
J_1&\leq& C\int^t_0|A_1r_n||(\psi_1(s,Y_n)-\psi_1(s,\check{Y}))h^1_n|ds\\
&\leq& C\int^t_0|A_1r_n|\|\psi_1(s,Y_n)-\psi_1(s,\check{Y})\|_{\mathcal{L}_2(U;V)}|h^1_n|_Uds\\
&\leq& C\int^t_0|A_1r_n|\|Y_n-\check{Y}\||h^1_n|_Uds\\
&\leq& C\int^t_0|A_1r_n|\|r_n+\eta_n\||h^1_n|_Uds\\
&\leq& \varepsilon\int^t_0\|r_n(s)\|^2_2ds+C\int^t_0(\|r_n\|^2+\|\eta_n\|^2)|h^1_n|^2_Uds,
\end{eqnarray*}
by \textbf{Hypothesis H0} and Corollary \ref{cor-3},
\begin{eqnarray*}
J_2&\leq& C\int^t_0\|\psi_1(s,\check{Y})(h^1_n-h_1)\|\|r_n\|ds\\
&\leq& C\int^t_0\|\psi_1(s,\check{Y})\|_{\mathcal{L}_2(U;V)}|h^1_n-h_1|_U\|r_n\|ds\\
&\leq& C\Big(\int^t_0\|r_n\|^2ds\Big)^{\frac{1}{2}}\Big(\int^t_0\|\psi_1(s,\check{Y})\|^2_{\mathcal{L}_2(U;V)}|h^1_n-h_1|^2_Uds\Big)^{\frac{1}{2}}\\
&\leq& C\Big(1+\sup_{t\in[0,T]}\|\check{Y}\|\Big)(2M)^{\frac{1}{2}}\Big(\int^t_0\|r_n\|^2ds\Big)^{\frac{1}{2}}.
\end{eqnarray*}
Collecting all estimates above, we have
\begin{eqnarray}\label{equ-80}
&&\|r_n(t)\|^2+\int^t_0 \|r_n(s)\|^2_2ds\\  \notag
&\leq& C\int^t_0(1+\|v_n\|^2)\|v_n\|^2_2\|r_n\|^2ds+C\int^t_0(1+\|\check{v}\|^2)\|\check{v}\|^2_2\|r_n\|^2ds\\ \notag
&& +C\int^t_0(|r_n|^2+|\nabla \eta_n|^2)ds+\int^t_0(\|r_n\|^2+\|\eta_n\|^2)|h^1_n|^2_Uds\\ \notag
&& +CM^{\frac{1}{2}}\Big(\int^T_0\|r_n\|^2ds\Big)^{\frac{1}{2}}.
\end{eqnarray}

\textbf{$H^1$ estimate of $\eta_n$}. \quad Similarly to the above, we omit the detail and only give the result,
\begin{eqnarray}\label{equ-81}
&&\|\eta_n(t)\|^2+\int^t_0 \|\eta_n(s)\|^2_2ds\\  \notag
&\leq& C\int^t_0(1+\|\check{T}\|\|\check{T}\|_2)\|r_n\|^2ds+C\int^t_0(1+\|v_n\|^2)\|v_n\|^2_2\|\eta_n\|^2ds\\ \notag
&& +C\int^t_0(\|r_n\|^2+\|\eta_n\|^2)|h^2_n|^2_Uds \notag
+CM^{\frac{1}{2}}\Big(\int^T_0\|\eta_n\|^2ds\Big)^{\frac{1}{2}}.
\end{eqnarray}
Thus, by (\ref{equ-80}) and (\ref{equ-81}),
\begin{eqnarray}\label{equ-82}
&&\|\rho_n(t)\|^2+\int^t_0 \|\rho_n(s)\|^2_2ds\\  \notag
&\leq&C\int^t_0(1+\|\check{v}\|^2)\|\check{v}\|^2_2\|r_n\|^2ds+ C\int^t_0(1+\|\check{T}\|\|\check{T}\|_2)\|r_n\|^2ds\\ \notag
&&+C\int^t_0(1+\|v_n\|^2)\|v_n\|^2_2\|\rho_n\|^2ds+C\int^t_0\|\rho_n\|^2|h_n|^2_Uds \\ \notag
&&+CM^{\frac{1}{2}}\Big[\Big(\int^T_0\|r_n\|^2ds\Big)^{\frac{1}{2}}+\Big(\int^T_0\|\eta_n\|^2ds\Big)^{\frac{1}{2}}\Big].
\end{eqnarray}
Applying Gronwall inequality to (\ref{equ-82}),
\begin{eqnarray*}
&&\sup_{t\in [0,T]}\|\rho_n(t)\|^2+\int^T_0 \|\rho_n(s)\|^2_2ds\\
&\leq& CM^{\frac{1}{2}}\Big[\Big(\int^T_0\|r_n\|^2ds\Big)^{\frac{1}{2}}+\Big(\int^T_0\|\eta_n\|^2ds\Big)^{\frac{1}{2}}\Big] \times \\ && \quad \exp\Big\{C\int^T_0
\Big[(1+\|\check{v}\|^2)\|\check{v}\|^2_2+(1+\|\check{T}\|\|\check{T}\|_2)+(1+\|v_n\|^2)\|v_n\|^2_2+|h_n|^2_U\Big]ds\Big\},
\end{eqnarray*}
moreover, Corollary \ref{cor-2} and Corollary \ref{cor-3} imply
\begin{eqnarray*}
\lim_{n\rightarrow \infty}\Big(\int^T_0\|r_n\|^2ds+\int^T_0\|\eta_n\|^2ds\Big)=0
\end{eqnarray*}
and
\begin{eqnarray*}
\exp\Big\{C\int^T_0
\Big[(1+\|\check{v}\|^2)\|\check{v}\|^2_2+(1+\|\check{T}\|\|\check{T}\|_2)+(1+\|v_n\|^2)\|v_n\|^2_2+|h_n|^2_U\Big]ds\Big\}
\leq C(T,\|Y_0\|, M),
\end{eqnarray*}
hence, we have
\begin{eqnarray}\label{equ-83}
|Y_n-\check{Y}|^2_{\Re}=\sup_{t\in [0,T]}\|\rho_n(t)\|^2+\int^T_0 \|\rho_n(s)\|^2_2ds\rightarrow 0, \quad n\rightarrow \infty.
\end{eqnarray}
$\hfill\blacksquare$
\section{Main Result}

\begin{thm}\label{thm-1}
Suppose that \textbf{Hypothesis H0} holds. Then for any $Y_0\in V$, $\{Y^\varepsilon\}$ satisfies the large deviation principle on $C([0,T],V)\cap L^2([0,T], D(A))$ with a good rate function given by (\ref{eq-5}).
\end{thm}

\textbf{Proof of Theorem \ref{thm-1} }
To prove the theorem, it suffices to verify the two conditions in \textbf{Hypothesis H1} so that Theorem \ref{thm-2} is applicable to obtain the large deviation principle for $Y^\varepsilon$.

\textbf{Step 1} \quad First, we show that the set $K_M=\{\mathcal{G}^0(\int^{\cdot}_{0}h(s)ds): h\in T_M\}$ is compact subset of $\Re$, where $\mathcal{G}^0$ is defined in (\ref{equ-28}).

Let $\{Y_{n}\}$ be a sequence in $K_M$ where $Y_{n}$ is the unique strong solution of (\ref{equ-9}) with $h_n\in T_M$. Keep in mind that we use the
weak topology on $T_M$. Hence there exists a subsequence (which we still denote it by $\{h_n\}$) converging to a limit $h$ weakly in $T_M$. Denote $Y_h$ be the strong solution of (\ref{equ-9}) with $h$. Corollary \ref{cor-1} and Theorem \ref{thm-6} establish that $Y_{n}\rightarrow Y_h$ in $\Re$ as $n\rightarrow \infty$, which implies that $K_M=\{\mathcal{G}^0(\int^{\cdot}_{0}h(s)ds): h\in T_M\}$ is compact subset of $\Re$.

%

\textbf{Step 2} \quad  Suppose that $\{h_{\varepsilon}: \varepsilon>0\}\subset \mathcal{A}_M$ for any fixed $M<\infty$ and $h_{\varepsilon}$ converge to $h$ as $T_M-$ valued random elements in distribution. Recall (\ref{eq def G epsilon}) the definition of $\mathcal{G}^{\varepsilon}$. Girsanov's theorem establishes that $\bar{Y}_{h_{\varepsilon}}=\mathcal{G}^{\varepsilon}(W(\cdot)+\frac{1}{\sqrt{\varepsilon}}\int^{\cdot}_0h^{\varepsilon}(s)ds)$
solves the following equation
\begin{eqnarray}\label{equ-19}
\left\{
  \begin{array}{ll}
    d\bar{Y}_{h_{\varepsilon}}(t)+A\bar{Y}_{h_{\varepsilon}}(t)dt+B(\bar{Y}_{h_{\varepsilon}}(t),\bar{Y}_{h_{\varepsilon}}(t))dt
    +G(\bar{Y}_{h_{\varepsilon}}(t))dt=\psi (\bar{Y}_{h_{\varepsilon}}) h_{\varepsilon}dt+\sqrt{\varepsilon}\psi(\bar{Y}_{h_{\varepsilon}}) dW(t), \\
    \bar{Y}_{h_{\varepsilon}}(0)=y_0.
  \end{array}
\right.
\end{eqnarray}
By $\rm It\hat{o}'s$ formula,
\begin{eqnarray}\label{eq P26 01}
\sup_{\epsilon\in(0,1)}\Big(\mathbb{E}\Big(\sup_{0\leq t\leq T}\|\bar{Y}_{h_{\varepsilon}}(t)\|^2+\int^T_0\|\bar{Y}_{h_{\varepsilon}}(t)\|^2_2dt\Big)\Big)\leq C<\infty.
\end{eqnarray}
Introducing an auxiliary process $Z_{\varepsilon}=(Z^1_{\varepsilon},Z^2_{\varepsilon})$,
\begin{eqnarray}\label{equ-20}
\left\{
  \begin{array}{ll}
    dZ_{\varepsilon}(t)+AZ_{\varepsilon}(t)dt=\sqrt{\varepsilon}\psi(t,\bar{Y}_{h_{\varepsilon}}) dW(t), \\
    Z_{\varepsilon}(0)=0.
  \end{array}
\right.
\end{eqnarray}
\textbf{Hypothesis H0} and (\ref{eq P26 01}) imply that
\begin{eqnarray}\label{eq P26 Z}
\lim_{\varepsilon\rightarrow 0}\mathbb{E}\Big(\sup_{0\leq t\leq T}\|Z_{\varepsilon}(t)\|^2+\int^T_0\|Z_{\varepsilon}(t)\|^2_2dt\Big)=0.
\end{eqnarray}
Since $T_M$ is a Polish space, by the Skorohod representation theorem, we can construct a stochastic basis $(\Omega^1,\mathcal{F}^1,\mathbb{P}^1)$ and, on this basis, $T_M\otimes T_M\otimes C([0,T],V)\cap L^2([0,T],D(A))$-valued random variables processes $(\tilde{h}_{\varepsilon},\tilde{h}, \tilde{Z}_{\varepsilon} )$ such that the joint distribution of
$(\tilde{h}_{\varepsilon},\tilde{Z}_{\varepsilon})$ is the same as that of $(h_{\varepsilon}, Z_{\varepsilon})$, $\tilde{Z}_{\varepsilon}\rightarrow 0$ a.s. in $C([0,T],V)\cap L^2([0,T],D(A))$, the distribution of $h$ coincides with that of $\tilde{h}$ and $\tilde{h}_{\varepsilon}\rightarrow \tilde{h}$ a.s. as $T_M-$ valued random elements.
Let ${X}_{\tilde{h}_{\varepsilon}}(t)$ be the solution of
\begin{eqnarray}\label{equ-21}
\left\{
  \begin{array}{ll}
    d{X}_{\tilde{h}_{\varepsilon}}(t)+A{X}_{\tilde{h}_{\varepsilon}}(t)dt
    +B({X}_{\tilde{h}_{\varepsilon}}(t)+\tilde{Z}_{\varepsilon},{X}_{\tilde{h}_{\varepsilon}}(t)+\tilde{Z}_{\varepsilon})dt
    +G({X}_{\tilde{h}_{\varepsilon}}(t)+\tilde{Z}_{\varepsilon})dt=\psi(t,{X}_{\tilde{h}_{\varepsilon}}+\tilde{Z}_{\varepsilon}) \tilde{h}_{\varepsilon}dt, \\
    {X}_{\tilde{h}_{\varepsilon}}(0)=y_0.
  \end{array}
\right.
\end{eqnarray}
The uniqueness of (\ref{equ-21}) implies that ${X}_{\tilde{h}_{\varepsilon}}$ has the same distribution with $\bar{Y}_{h_{\varepsilon}}-Z_{\varepsilon}$.
Using similar arguments as in Sect. \ref{Sec 4}, we can prove
\begin{eqnarray*}
{X}_{\tilde{h}_{\varepsilon}}\rightarrow {X}_{\tilde{h}}\text{ in }\Re,\ \ \ \mathbb{P}^1\text{-a.s.}
\end{eqnarray*}
which satisfies
\begin{eqnarray*}
\left\{
  \begin{array}{ll}
    d{X}_{\tilde{h} }(t)+A{X}_{\tilde{h} }(t)dt
    +B({X}_{\tilde{h} }(t) ,{X}_{\tilde{h} }(t) )dt
    +G({X}_{\tilde{h} }(t) )dt=\psi(t,{X}_{\tilde{h} }) \tilde{h} dt, \\
    {X}_{\tilde{h} }(0)=y_0.
  \end{array}
\right.
\end{eqnarray*}
Recall (\ref{equ-28}) the definition of $\mathcal{G}^0$. Combining ${X}_{\tilde{h}_{\varepsilon}}$ has the same distribution with $\bar{Y}_{h_{\varepsilon}}-Z_{\varepsilon}$ and (\ref{eq P26 Z}), we obtain \textbf{Hypothesis H1} (i).
The proof is complete.
$\hfill\blacksquare$


\def\refname{ References}

\end{document}